
\documentstyle{amsppt}
\magnification=\magstep1
\NoRunningHeads
\pageno=1
\vsize=7.4in


\def\bull{\vrule height .9ex width .8ex depth -.1ex}

\def\Ave{\mathop{\text {Ave}}}
\def\Pr{\text{Pr}}

\topmatter
\title
Subspaces of rearrangement-invariant spaces
\endtitle

\author
Francisco L. Hernandez and
Nigel J. Kalton
\endauthor

\address1 Facultad de Matematicas, Universidad Complutense,
28040 Madrid, Spain \endaddress \email pacoh\@mat.ucm.es
\endemail

\address2 Department of Mathematics, University of Missouri
Columbia, MO 65211, U.S.A.  \endaddress \email
mathnjk\@mizzou1.bitnet \endemail \thanks The first author
was partially supported by DGICYT grant PB-910377.  The
second author was supported by NSF grant DMS-9201357.
\endthanks \subjclass 46B03 \endsubjclass \abstract We prove
a number of results concerning the embedding of a Banach
lattice $X$ into an r.i. space $Y$.  For example we show
that if $Y$ is an r.i. space on $[0,\infty)$ which is
$p$-convex for some $p>2$ and has nontrivial concavity then
any Banach lattice $X$ which is $r$-convex for some $r>2$
and embeds into $Y$ must embed as a sublattice.  Similar
conclusions can be drawn under a variety of hypotheses on
$Y$; if $X$ is an r.i. space on $[0,1]$ one can replace the
hypotheses of $r$-convexity for some $r>2$ by $X\neq L_2.$

We also show that if $Y$ is an order-continuous Banach
lattice which contains no complemented sublattice
lattice-isomorphic to $\ell_2,$ $X$ is an order-continuous
Banach lattice so that $\ell_2$ is not complementably
lattice finitely representable in $X$ and $X$ is isomorphic
to a complemented subpace of $Y$ then $X$ is isomorphic to a
complemented sublattice of $Y^N$ for some integer $N.$

\endabstract

\endtopmatter

\document \baselineskip=14pt

\heading {1.  Introduction}\endheading \vskip10pt

The study of the Banach space geometry of general
rearrangement-invariant Banach function spaces may be
considered to originate with the work of Bretagnolle and
Dacunha-Castelle on subspaces of Orlicz function spaces
\cite{3}.  A very important development in the theory was
the publication of a systematic study of r.i. spaces by
Johnson, Maurey, Schechtman and Tzafriri in 1979 \cite{21}.
The appearance of this memoir revolutionized the subject.
Since then, a number of authors have considered problems of
classifying subspaces of certain special r.i. spaces; see
\cite{5}, \cite{6}, \cite{7} ,\cite{8}, \cite{9}, \cite{13},
\cite{14}, \cite{17}, \cite {19}, \cite{20}, \cite{39},
\cite{40} for a variety of different results of this type.

In general, most of the literature relates to the problem of
embedding a Banach lattice $X$ (either atomic or nonatomic)
with additional symmetry conditions into an r.i. space $Y,$
and the techniques used rely heavily on symmetrization.  In
\cite{27}, however, the second author considered the general
problem of determining conditions when an order-continuous
Banach lattice $X$ could be complementably embedded in an
order-continuous Banach lattice $Y,$ minimizing the use of
symmetry.  The aim was to show that under certain hypotheses
on $X$ and $Y$ one could deduce that $X$ (or perhaps only a
non-trivial band in $X$) would be lattice-isomorphic to a
complemented sublattice of $Y.$ A number of such results
were obtained (we refer for details to \cite{27}); of
course, the additional assumption that either $X$ or $Y$ is
r.i. could still be used to obtain stronger results of this
nature.  In the final section of this paper (Section 8,
which can be read independently of the remainder) we obtain
a significant improvement of one of the results of \cite{27}
by showing that if $X,Y$ are order-continuous separable
Banach lattices, such that $Y$ contains no complemented
sublattice which is lattice-isomorphic to $\ell_2$ and
$\ell_2$ is not complementably lattice finitely
representable in $X$, and if $X$ is isomorphic to a
complemented subspace of $Y$ then $X$ is lattice-isomorphic
to a complemented sublattice of $Y^N$ for some $N.$ Of
course if $Y$ is r.i. then $X$ must be a complemented
sublattice of $Y$ itself.

The main body of the paper (Sections 3-7) is concerned with
similar problems but without assumptions of complementation.
We consider an r.i. space $Y$ on $[0,\infty)$ (or $[0,1]$,
but there our results are not quite so strong) and consider
a generally nonatomic Banach lattice $X$ which is isomorphic
to a subspace of $Y$; we would like to show, under
appropriate hypotheses that $X$ is lattice-isomorphic to a
sublattice of $Y.$ Of course, there is no hope of such a
result in general; the spaces $L_p[0,1]$ for $1\le p<2$ have
a very rich subspace structure (cf.  \cite{39}, \cite{40});
in particular $L_r$ embeds into $L_p$ if $p<r\le 2.$
However, there are some suggestive results in the literature
which tend to indicate the possibility of strong conclusions
if $Y$ is ``on the other side of 2.''

We first observe that Johnson, Maurey, Schechtman and
Tzafriri \cite{21} Theorem 1.8, showed that if $X$ is a
Banach lattice which embeds into $L_p[0,1]$ where $p>2$ and
$X$ is $r$-convex for some $r>2$ (or, equivalently $\ell_2$
is not lattice finitely representable in $X$) then $X$ is
lattice-isomorphic to $L_p(\mu)$ for some measure $\mu,$ and
so is lattice-isomorphic to a sublattice of $L_p.$ Note that
this result requires no symmetry conditions on $X.$ For the
case when $X$ is an r.i. space on $[0,1]$ there are some
other positive results.  In \cite{21} Theorem 7.7 shows that
if $Y=L_F[0,\infty)$ is a $p$-convex Orlicz space, with
nontrivial concavity, where $p>2$ and if $X$ is an r.i.
space on $[0,1]$ which embeds into $Y$, with $X\neq
L_2[0,1],$ then $X$ must be lattice-isomorphic to a
sublattice of $Y.$ Later Carothers \cite{5} proved the same
result for the Lorentz spaces $L_{p,q}$ where $2<q<p.$ These
spaces are also strictly 2-convex (i.e.  $r$-convex for some
$r>2$).  However in \cite{6}, Carothers extended his work to
the Lorentz spaces $L_{p,q}$ where $1\le q\le 2<p.$ These
spaces are not even 2-convex.

Our main results include all these previous theorems.  In
Theorem 7.2, we show that if $Y$ is a strictly 2-convex r.i.
space on $[0,\infty)$ with nontrivial concavity and $X$ is a
strictly 2-convex Banach lattice then if $X$ embeds into
$Y,$ then $X$ is lattice-isomorphic to a sublattice of $Y$.
The assumption of strict 2-convexity on $Y$ can be relaxed
for a special class of r.i. spaces which we term of
Orlicz-Lorentz type (this class includes all reflexive
Orlicz and Lorentz spaces); if $Y$ is of Orlicz-Lorentz type
we need only assume that $Y$ is 2-convex or that its lower
Boyd index $p_Y>2.$ In the case when $Y$ is an r.i. space on
$[0,1]$ our results are not quite as good; for example if
$Y$ is strictly 2-convex and has nontrivial concavity and
$X$ is strictly 2-convex we deduce only that some nontrivial
band in $X$ is lattice-isomorphic to a sublattice of $Y$.
In the case when $X$ is an r.i. space on $[0,1]$ we give
(Corollary 7.4) a very general result which includes the
above mentioned results of \cite {5}, \cite{6} and \cite{21}
for Orlicz and Lorentz spaces.  Precisely, suppose $Y$ is an
r.i. space on $[0,1]$ or $[0,\infty)$ with nontrivial
concavity and suppose that either $Y$ is strictly 2-convex
or $Y$ is of Orlicz-Lorentz type with $p_Y>2$; suppose $X$
is an r.i. space on $[0,1]$ which embeds into $Y$.  Then
either $X=L_2[0,1]$ or $X$ is lattice-isomorphic to a
sublattice of $Y$ (so that $X=Y_f[0,1],$ for some $f\in Y).$

We also give a result on embedding $L_p[0,1]$ where $p>2$
into a $p$-concave r.i. space $Y.$ We show in Theorem 7.7
that this implies that either the Haar basis of $L_p$ is
lattice finitely representable in $Y$ or $Y[0,1]=L_p[0,1].$
The former alternative is impossible if $Y$ is of
Orlicz-Lorentz type or is strictly 2-convex.

Let us now briefly discuss the method of proof of these
results.  For reasons discussed below, we consider
quasi-Banach lattices and develop a theory of
cone-embeddings.  If $X$ and $Y$ are quasi-Banach lattices,
a cone-embedding $L:X\to Y$ is a positive linear operator
such that for some $\delta>0,$ $\|Lx\|_Y \ge \delta\|x\|_X$
for every $x\ge 0.$ We consider cone-embeddings in Sections
4 and 5. The aim is to produce conditions on $X$ and $Y$ so
that one can pass from the existence of a cone-embedding to
the existence of a lattice-embedding.  Crucial use is made
of the theory of random measure representations of positive
operators.  A typical result is that if $X$ is strictly
1-convex and if $Y$ is an r.i. space on $[0,\infty)$ which
is an interpolation space between $L_1$ and $L_{\infty}$
then if $X$ cone-embeds into $Y$ it also lattice-embeds.
The assumption on $Y$ is satisfied if $Y$ is a Banach r.i.
space, by the Calder\'on-Mityagin theorem, but also holds
for certain non-Banach examples, where the lower Boyd index
$p_Y>1.$

The next step carried out in Section 6 is to consider the
case when $X$ is a Banach lattice which embeds into an r.i.
space $Y$.  The aim here is to put hypotheses on $X$ and $Y$
so that one can induce a cone-embedding $L:X_{1/2}\to
Y_{1/2}$ where $X_{1/2},Y_{1/2}$ are the 2-concavifications
of $X$ and $Y$ (these spaces may not be locally convex).
This can be done if one puts a somewhat technical hypothesis
on $X$ and $Y$ (Theorem 6.7).  To put this hypothesis in
perspective, let us note that if $X$ is an r.i. space on
$[0,1]$ and one aimed simply to guarantee that $L\neq 0$ it
would suffice to assume that the Haar basis of $X$ was not
equivalent to a disjoint sequence in $Y$.  This is a typical
hypothesis in \cite{21} (Theorems 5.1 and 6.1) where the aim
is only to draw the weaker conclusion that $X[0,1]\subset
Y[0,1].$ In fact some (and perhaps all) of these results can
be recovered from our method.  However, to obtain $X$ as a
sublattice we need $L$ to be a cone-embedding.  Fortunately
our stronger technical condition is satisfied when $Y$ is
strictly 2-convex or of Orlicz-Lorentz type.

Finally one can put these steps together and obtain, under
the right hypotheses, that if $X$ embeds into $Y$ then
$X_{1/2}$ lattice-embeds into $Y_{1/2}$ and so $X$
lattice-embeds into $Y.$

This research was carried out during a visit of the first author to the
University of Missouri in October 1993 and a visit of the second author
to the Complutense University in Madrid in June 1994.

\vskip10pt \heading{2.  Definitions and notation}\endheading
\vskip10pt

We first recall that a (quasi-)Banach lattice $X$ is said to
be order-continuous if and only if every order-bounded
increasing sequence is norm convergent (see \cite{34} p.7).
A quasi-Banach lattice which does not contain a copy of
$c_0$ is automatically order-continuous but the converse is
false.  An atom in a Banach lattice is a positive element
$a$ so that $0\le x\le a$ implies that $x=\alpha a$ for some
$0\le \alpha\le 1.$ A Banach lattice is nonatomic if it
contains no atoms.  The reader is referred to
Lindenstrauss-Tzafriri \cite{34} or Meyer-Nieberg \cite{36}
as a general reference for Banach lattices.

We will in general use the same notation as in \cite{27}.
Let $\Omega$ be a Polish space (i.e. a separable complete
metric space) and let $\mu$ be a $\sigma-$finite Borel
measure on $\Omega.$ We refer to the pair $(\Omega,\mu)$ as
a Polish measure space; if $\mu$ is a probability measure
then we say $(\Omega,\mu)$ is a Polish probability space.
If $E$ is a Borel set then $\chi_E$ denotes its indicator
function.  We denote by $L_0(\mu)$ the space of all Borel
measurable functions on $\Omega$, where we identify
functions differing only on a set of measure zero; the
natural topology of $L_0$ is convergence in measure on sets
of finite measure.  If $0<p\le 1,$ an admissible $p$-norm is
then a lower-semi-continuous map $f\to \|f\|$ from
$L_0(\mu)$ to $[0,\infty]$ such that:\newline \vskip.2truecm
\flushpar (a) $\|\alpha f\| =|\alpha|\|f\|$ whenever
$\alpha\in\bold R,\ f\in L_0.$ \newline (b) $\|f+g\|^p\le
\|f\|^p+\|g\|^p,$ for $f,g\in L_0.$ \newline (c) $\|f\| \le
\|g\|,$ whenever $|f| \le |g|$ a.e.  (almost
everywhere).\newline (d) $\|f\|<\infty$ for a dense set of
$f\in L_0,$ \newline (e) $\|f\|=0$ if and only if $f=0$ a.e.
\newline \vskip.1truecm If $p=1,$ we call $\|\,\|$ an
admissible norm; an admissible quasinorm is an admissible
$p$-norm for some $0<p\le 1.$

A quasi-K\"othe function space on $(\Omega,\mu)$ is defined
to be a dense order-ideal $X$ in $L_0(\mu)$ with an
associated admissible quasinorm $\|\,\|_X$ such that if
$X_{max}=\{f:\|f\|_X<\infty\}$ then either:  \newline (1)
$X=X_{max}$ ($X$ is {\it maximal}) or:\newline (2) $X$ is
the closure of the simple functions in $X_{max}$ ($X$ is
{\it minimal}).\newline

\vskip.1truecm If $\|\,\|_X$ is a norm then $X$ is called a
K\"othe function space.  Notice that according to our
description we consider $\|\,\|_X$ to be well-defined on
$L_0.$ Any order-continuous K\"othe function space is
minimal.  Also any K\"othe function space which does not
contain a copy of $c_0$ is both maximal and minimal.

Given any K\"othe function space $X$ and $0<p<\infty$ we
define $X_p$ to be the quasi-K\"othe space of all $f$ such
that $|f|^{p}\in X$ with the associated admissible quasinorm
$\|f\|_{X_p}=\||f|^p\|_X^{1/p}.$ It is readily verified that
$\|\,\|_{X_p}$ is an admissible $p$-norm when $0<p<1$ and an
admissible norm when $p>1.$ We will primarily use the case
$p=1/2$ in this paper.  We will also use the subscript $+$
to denote the positive cone in a variety of situations, e.g.
$X_+=\{f:f\in X,\ f\ge 0\}.$

If $X$ is an order-continuous K\"othe function space then
$X^*$ can be identified with the K\"othe function space of
all $f$ such that:  $$ \|f\|_{X^*}=\sup_{\|g\|_X\le 1} \int
|fg|\,d\mu <\infty.$$ $X^*$ is always maximal.

If $\mu$ is a probability measure then we say following
\cite{21}, that a K\"othe function space $X$ is {\it good\/}
if $L_{\infty}\subset X\subset L_1$ and further for $f\in
L_0,$ $\|f\|_1\le \|f\|_X\le 2\|f\|_{\infty}.$ It is
well-known that any separable order-continuous Banach
lattice can be represented as (i.e. is isometrically
lattice-isomorphic to) a good K\"othe function space on some
Polish probability space $(\Omega,\mu)$ (see \cite{21} and
\cite{34}).

In the case when $X$ is nonatomic we can require that
$\Omega=[0,1]$ and $\mu=\lambda$ is Lebesgue measure.
Alternatively we can take $\Omega=\Delta=\{-1,+1\}^{\bold
N}$ to be the Cantor group and take $\mu$ to be normalized
Haar measure on $\Delta$ which we again denote by $\lambda.$
We will use this second representation freely and now take
the opportunity to introduce some notation from \cite{27}.

Thus for $\epsilon_k=\pm1,$ we denote by
$\Delta(\epsilon_1,\ldots,\epsilon_n)$ the clopen subset of
$\Delta$ of all $(d_j)_{j=1}^{\infty}$ such that
$d_j=\epsilon_j$ for $1\le j\le n.$ For each $n$ let $\Cal
A_n$ denote the collection of
$\Delta(\epsilon_1,\ldots,\epsilon_n)$.  Let $CS_n$ denote
the linear span of $\{\chi_E:E\in\Cal A_n\}.$ We also define
the Haar functions
$h_E=\chi_{\Delta(\epsilon_1,\ldots,\epsilon_n,+1)}-
\chi_{\Delta(\epsilon_1,\ldots,\epsilon_n,-1)}$ for
$E=\Delta(\epsilon_1,\ldots,\epsilon_n).$

A K\"othe function space (or, more generally a quasi-K\"othe
function space) $X$ is said to be $p-$convex (where
$0<p<\infty)$ if there is a constant $C$ such that for any
$f_1,\ldots,f_n\in X$ we have $$
\|(\sum_{i=1}^n|f_i|^p)^{1/p}\|_X \le
C(\sum_{i=1}^n\|f_i\|_X^p)^{1/p}.$$ $X$ is said to have an
upper $p$-estimate if for some $C$ and any disjoint
$f_1,\ldots,f_n\in X,$ $$ \|\sum_{i=1}^n f_i\|_X \le
C(\sum_{i=1}^n\|f_i\|_X^p)^{1/p}.$$ $X$ is said to be
$q-$concave ($0<q<\infty$) if for some $c>0$ and any
$f_1,\ldots,f_n\in X$ we have $$
\|(\sum_{i=1}^n|f_i|^q)^{1/q}\|_X \ge
c(\sum_{i=1}^n\|f_i\|_X^q)^{1/q}.$$ $X$ is said to have a
lower $q$-estimate if for some $c>0$ and any disjoint
$f_1,\ldots,f_n\in X,$ $$ \|\sum_{i=1}^n f_i\|_X \ge
c(\sum_{i=1}^n \|f_i\|_X^q)^{1/q}.$$ Notice that a
quasi-K\"othe function space which satisfies a lower
$q$-estimate is automatically both maximal and minimal since
it cannot contain a copy of $c_0.$ A K\"othe function space
must, of course, be 1-convex.  A quasi-K\"othe function
space must satisfy an upper $p$-estimate for some $p>0$ but
need not be $p$-convex for any $p>0;$ however, if $X$
satisfies a lower $q$-estimate for some $q<\infty$ then it
is $p$-convex for some $p>0.$ This result is proved in
\cite{24} (Theorems 4.1 and 2.2) and a simpler proof is
presented in \cite{30} Theorem 3.2.  A quasi-K\"othe
function space which is $s$-convex for some $s>0$ and
satisfies an upper $r$-estimate is $p$-convex for every
$0<p<r$ (see \cite{24}).

A (quasi-)Banach lattice $X$ is $p$-convex, satisfies an
upper $p$-estimate, is $q$-concave or satisfies a lower
$q$-estimate according as any concrete representation of $X$
as a K\"othe function space has the same property.  We shall
say that $X$ is {\it strictly p-convex} if it is $r-$convex
for some $r>p$ and {\it strictly q-concave} if it is
$s$-concave for some $s<q.$

A Banach space $X$ is said to be of (Rademacher) type $p$
($1\le p\le 2$) if there is a constant $C$ so that for any
$x_1,\ldots,x_n \in X,$ $$ \mathop{\text{Ave
}}_{\epsilon_i=\pm 1}\|\sum_{i=1}^n \epsilon_i x_i\| \le
C(\sum_{i=1}^n\|x_i\|^p)^{1/p}$$ and $X$ is of cotype $q$
($2\le q<\infty$) if for some $c>0$ and any
$x_1,\ldots,x_n\in X$ we have $$ \mathop{\text{Ave
}}_{\epsilon_i=\pm1}\|\sum_{i=1}^n\epsilon_ix_i\| \ge
c(\sum_{i=1}^n\|x_i\|^q)^{1/q}.$$ We recall that a
(quasi-)Banach lattice has nontrivial cotype (i.e. has
cotype $q<\infty$ for some $q$) if and only if it has
nontrivial concavity (i.e. is $q-$concave for some
$q<\infty$).  If $X$ is a Banach lattice which has
nontrivial concavity then there is a constant $C=C(X)$ so
that for any $x_1,\ldots,x_n\in X$ we have $$ \frac1C
(\Ave_{\epsilon_k=\pm1}\|\sum_{k=1}^n\epsilon_kx_k\|^2)^{1/2}\le
\|(\sum_{k=1}^n|x_k|^2)^{1/2}\|_X \le C
(\Ave_{\epsilon_k=\pm1}\|\sum_{k=1}^n\epsilon_kx_k\|^2)^{1/2}.$$
In fact we will need the same conclusion for quasi-Banach
lattices; as far as we know this has never been explicitly
stated although it is probably well-known.  We therefore
state it formally as a Proposition.

\proclaim{Proposition 2.1}Let $X$ be a quasi-Banach lattice
with nontrivial concavity (equivalently nontrivial cotype).
Then there is a constant $C=C(X)$ so that for any
$x_1,\ldots,x_n\in X$ we have $$ \frac1C
(\Ave_{\epsilon_k=\pm1}\|\sum_{k=1}^n\epsilon_kx_k\|^2)^{1/2}\le
\|(\sum_{k=1}^n|x_k|^2)^{1/2}\|_X \le C
(\Ave_{\epsilon_k=\pm1}\|\sum_{k=1}^n\epsilon_kx_k\|^2)^{1/2}.$$
\endproclaim

\demo{Proof}We have that $X$ is $q$-concave for some
$q<\infty$.  As remarked above it is also $p$-convex for
some $p>0.$ It is now easy to adapt the standard argument
based on Khintchine's inequality as in \cite{34} Theorem
1.d.6, p. 49.\bull\enddemo

\demo{Remark}In fact we will only apply this Proposition in
situations when the $p$-convexity of $X$ for some $p>0$ is
automatic (i.e.  $X$ is the concavification of some K\"othe
function space).\enddemo

Let us now turn to rearrangement-invariant spaces (cf.
\cite{21},\cite{34}).  For any $f\in L_0(\Omega,\mu)$ we
define its decreasing rearrangement $f^*\in
L_0[0,\mu(\Omega))$ by $f^*(t)=\inf\{x:\mu(|f|>x)\le t\}.$
Now let $X$ be a quasi-K\"othe function space on either
$[0,\infty)$ or $[0,1]$ with Lebesgue measure.  We say that
$X$ is a quasi-Banach {\it rearrangement-invariant} (r.i.)
space if $\|f\|_X=\|f^*\|_X$ for all $f\in L_0,$ and if
$\|\chi_{[0,1]}\|_X=1.$ We use the term r.i. space for a
Banach r.i. space.  If $X$ is a quasi-Banach r.i. space on
$[0,\infty)$ (respectively, $[0,1]$) and $(\Omega,\mu)$ is a
Polish measure space (respectively, with $\mu(\Omega)\le
1,$) then we define $X(\Omega,\mu)$ to be the set of $f\in
L_0(\mu)$ such that $f^*\in X$ with $\|f\|_X= \|f^*\|_X.$
For example, it will be of some advantage to consider
$X(\Delta,\lambda)$ in place of $X[0,1].$ Let us remark that
if $X$ is a quasi-Banach r.i. space on $[0,1]$ then it is
always possible to write $X=Y[0,1]$ where $Y$ is some
quasi-Banach r.i. space on $[0,\infty).$ We will only be
interested in separable (or order-continuous) r.i. spaces,
which are necessarily minimal.

On any quasi-Banach r.i. space $X$ on $[0,\infty)$ (resp.
$[0,1]$) we define the dilation operators $D_s$ for
$0<s<\infty$ by $$ D_sf(t) = f(t/s)$$ for all $t$ (resp.
whenever $0\le t\le \min(1,s)$ and $D_sf(t)=0$ otherwise).
The Boyd indices $p_X$ and $q_X$ are defined by $$ \align
p_X &=\lim_{s\to\infty}\frac{\log s}{\log \|D_s\|}\\ q_X &=
\lim_{s\to 0}\frac{\log s}{\log \|D_s\|}.  \endalign $$ In
general $0< p_X \le q_X\le\infty;$ if $X$ is a Banach r.i.
space (i.e. is 1-convex) then $1\le p_X.$ If $X$ is an
order-continuous Banach r.i. space, then $X$ has an
unconditional basis if and only if $1<p_X\le q_X<\infty;$ in
this case the Haar basis of $X$ is an unconditional basis
(see \cite{34} p. 157-161).

Recall that if $f\in L_0(\Omega,\mu)$ then
$f^{**}(t)=\frac1t\int_0^tf^*(s)ds,$ for $t>0.$ We say that
a quasi-Banach r.i. space $X$ on $[0,1]$ or $[0,\infty)$ has
property (d) if there exists $C$ so that if $f\in X$ and
$g\in L_0$ satisfy $g^{**}\le f^{**}$ then $g\in X$ with
$\|g\|_X\le C\|f\|_X.$ It is well-known that every Banach
r.i. space satisfies property (d) (cf.  \cite{34} p.125)
with $C=1.$ However there are non-locally convex examples;
any quasi-Banach r.i. space $X$ with $p_X>1$ satisfies
property (d) (see \cite{26}).  A quasi-Banach r.i. space
with property (d) is an interpolation space for the pair
$(L_1,L_{\infty});$ this is a mild generalization of the
classical Calder\'on-Mityagin theorem (\cite{4}, \cite{35})
which follows from considerations of the K-functional (see,
for example Bennett-Sharpley \cite{2}, Chapters 3 and 5;
this treats only the normed case, but the modifications are
trivial).

We also recall a definition from \cite{29}.  If $X$ is an
r.i. space on $[0,\infty)$ (resp.  $[0,1]$) we define $E_X$
to be the closed subspace of $X$ spanned by the functions
$e_n=\chi_{[2^{n},2^{n+1})}$ for $n\in\bold Z$ (resp.  $n\in
\bold Z_-=\{n:n<0\}$).  If $X$ is separable then $(e_n)$
forms an unconditional basis for $E_X$ and $E_X$ can be
regarded as a sequence space modelled on $\bold J=\bold Z$
or $\bold Z_-.$ We shall say that $X$ is of {\it
Orlicz-Lorentz type} if $E_X$ is naturally isomorphic to a
modular sequence space, i.e. there exist Orlicz functions
$(F_n)_{n\in\bold J}$ so that $E_X=\ell_{(F_n)} (\bold J)$
(see \cite{33} pp. 168ff).  This is a convenient definition
to specify a class of spaces $X$ which includes the standard
Orlicz spaces and Lorentz spaces, and a variety of ``mixed''
spaces.

To illustrate these ideas consider the following method of
defining an r.i. space on $[0,\infty).$ Let $Y$ be a K\"othe
function space on $[0,\infty)$ with the property that the
dilation operators $D_t:Y\to Y$ are all bounded.  Then we
can define $p_Y,q_Y$ as in the rearrangement-invariant case.
Assume that $1<p_Y\le q_Y<\infty.$ Now let $\tilde Y$ be the
space defined by $f\in\tilde Y$ if and only if $f^*\in Y$
and define $\|f\|_{\tilde Y}=\|f^*\|_Y$.  The inequality
$(f+g)^* \le 2D_2f^*+2D_2g^*$ shows that $\|\,\|_{\tilde Y}$
is a quasinorm and that $\tilde Y$ is an order-ideal.  In
fact, we also have:

\proclaim{Proposition 2.2}There exists a constant $C$ so
that if $f\in L_0$ then $$ \|f\|_{\tilde Y}\le
\|\sum_{n\in\bold Z}f^{**}(2^{n})e_n\|_Y \le C\|f\|_{\tilde
Y}.$$\endproclaim

\demo{Proof}(Due to S. Montgomery-Smith).  Clearly $f^*\le
\sum_{n\in\bold Z}f^{**}(2^{n})e_n.$ However $f^{**}(2^n)
\le \sum_{k=1}^{\infty}2^{-k}f^*(2^{n-k}).$ Hence
$\sum_{n\in\bold Z}f^{**}(2^n)e_n \le
\sum_{k=1}^{\infty}2^{-k}D_{2^{k+1}}f^*$.  But now since
$p_Y>1$ it follows that
$\sum_{k=1}^{\infty}2^{-k}\|D_{2^{k+1}}\|_Y <\infty$ and the
result follows.\bull\enddemo

The proof above only uses the hypothesis that $p_Y>1,$ and
not that $q_Y<\infty.$ Proposition 2.2 shows that $\tilde Y$
is a Banach r.i. space by providing an equivalent norm.  It
is now immediate that $p_Y\le p_{\tilde Y}.$ We next show
that $E_{\tilde Y}$ coincides with $E_Y.$ This implies that
if $Y$ is an Orlicz-Musielak space or generalized Orlicz
space (cf.  \cite{37}) then the associated r.i. space
$\tilde Y$ is of Orlicz-Lorentz type as defined above.  In
particular, if we take $Y$ to be a weighted $L_p-$space
(with, of course the conditions $1<p_Y\le q_Y<\infty$) we
obtain the usual Lorentz spaces as examples of spaces of
Orlicz-Lorentz type.

\proclaim{Proposition 2.3}We have $E_{\tilde Y}=E_Y$ (and
the norms are equivalent).\endproclaim

\demo{Proof}In fact suppose $f=\sum_{n\in\bold Z}a_ne_n$
where $a_n\ge 0$ is finitely nonzero.  Let
$g=\sum_{n=0}^{\infty}D_{2^{-n}}f.$ The assumption
$q_Y<\infty$ and the fact that $q_{\tilde Y}\le q_Y$ is
sufficient to establish that $\|g\|_Y\le C\|f\|_Y$ and
$\|g\|_{\tilde Y}\le C\|f\|_{\tilde Y}$ for a suitable
constant $C.$ Note that $\|g\|_Y=\|g\|_{\tilde Y}$ since $g$
is decreasing.  The result follows immediately.\bull\enddemo

\vskip10pt

\heading{3.  Remarks on sublattices}\endheading \vskip10pt

In this section, we collect together some elementary remarks
on the structure of sublattices of r.i. spaces.

\proclaim{Lemma 3.1}Suppose $X$ is a quasi-K\"othe function
space on $(\Omega,\mu)$ and that $Y$ is a quasi-Banach r.i.
space on $[0,\infty).$ Suppose $I=[0,1]$ or $[0,\infty)$ and
that $U:X\to Y(I)$ is a lattice homomorphism.  Then there is
a lattice homomorphism $V:X\to Y(\Omega\times [0,\infty))$
so that for any $x\in X$ and $\alpha>0$ we have
$$\frac12\lambda(U|x|>2\alpha) \le
(\mu\times\lambda)(V|x|>\alpha)\le \lambda(U|x|>\alpha),$$
and such that $V$ can be represented as
$Vx(\omega,t)=a(\omega,t)x(\omega)$ where $a$ is a
nonnegative Borel function on $\Omega\times [0,\infty)$ of
the form $$a(\omega,t)=\sum_{k\in\bold Z}
2^{m(k,\omega)}e_k(t)$$ with $m:\bold Z\times \Omega\to
\bold Z\cup{-\infty}$ is a Borel map with $k\to m(k,\omega)$
decreasing for each $\omega.$ Furthermore if $I=[0,1]$ then
$a$ is supported on a set of measure one in the product
space.\endproclaim

\demo{Proof}It will suffice to consider the case when $X$
contains $L_{\infty}.$ We suppose the existence of a lattice
embedding $Ux=bx\circ\sigma$ where $b$ is a nonnegative
Borel function and $\sigma:I\to \Omega$ is a Borel map.
First pick $b'$ with $\frac12 b\le b'\le b$ so that
$b'=\sum_{n\in\bold Z}2^n\chi_{E_n}$ where $E_n$ are
disjoint Borel sets.  Let $U'x=b'x\circ\sigma.$

Now for each $n$ define the measure
$\nu_n(B)=\lambda(\cup_{k\ge n}E_k\cap \sigma^{-1}B).$ Since
$U'\chi_{\Omega}\in Y$ it is clear that each $\nu_n$ is a
finite measure.  Furthermore, if $\mu B=0$ then $U\chi_B=0$
a.e. and hence $\nu_n(B)=0.$ Hence we can find nonnegative
Borel functions $w_n$ on $\Omega$ so that $\nu_n(B)=\int_B
w_n d\mu,$ and we may suppose that $w_n(\omega)$ is
decreasing for each fixed $\omega.$ Notice that
$\int_{\Omega}w_nd\mu=\nu_n(\Omega)\le \lambda(I),$ so that
if $I=[0,1]$ then $\int_{\Omega}w_nd\mu\le 1$ for all $n,$

For any fixed $n\in\bold Z,$ we define $A_n=\{(\omega,t):\
t\le w_n(\omega)\}$ and let $a'=\sum_{n\in\bold
Z}2^n(\chi_{A_n}-\chi_{A_{n+1}}).$ Define $V':X\to
Y(\Omega\times (0,\infty))$ by
$V'x(\omega,t)=a'(\omega,t)x(t).$ Finally define $a$ a Borel
function on $\Omega\times (0,\infty)$ by setting
$a(\omega,t)=2^m$ if $a'(\omega,2^{k+1})=2^m$ where
$2^{k}\le t< 2^{k+1}$ and $k,m\in\bold Z.$ We set
$a(\omega,t)=0$ if $a'(\omega,2^{k+1})=0.$ Notice that
$(\mu\times\lambda)\{a>0\}\le (\mu\times\lambda)\{a'>0\} \le
1$ if $I=[0,1].$ Define $Vx(\omega,t)=a(\omega,t)x(t).$

Now suppose $x\ge 0,x\in X.$ Then $0\le Vx\le V'x.$
Furthermore for fixed $\omega,$ $$\lambda\{t:\
V'x(\omega,t)>\alpha\}\le 2\lambda\{t:\
Vx(\omega,t)>\alpha\}$$ so that
$$(\mu\times\lambda)(Vx>\alpha) \ge
\frac12(\mu\times\lambda)(V'x>\alpha).$$ Now again for fixed
$\alpha,$  let $F_n=\{2^nx<\alpha\le 2^{n+1}x\}.$ We
note that $$ \align (\mu\times\lambda)(V'x>\alpha)&=
\sum_{n\in\bold Z}\int_{F_n} w_nd\mu \\ &= \sum_{n\in\bold
Z}\nu_n(F_n)\\ &= \sum_{n\in\bold Z}\lambda(\bigcup_{k\ge
n}E_k\cap\sigma_n^{-1}F_n)\\ &= \lambda(U'x>\alpha).
\endalign $$ Hence
$$\frac12\lambda(U'x>\alpha)\le(\mu\times\lambda)(Vx>\alpha)
\le \lambda(U'x>\alpha).$$ Since $\frac12Ux\le U'x\le Ux$
the result follows.\bull \enddemo

We next state the immediate conclusion for lattice
embeddings.

\proclaim{Proposition 3.2}Let $X$ be a quasi-K\"othe
function space on $(\Omega,\mu).$ Suppose $Y$ is a
quasi-Banach r.i. space on $[0,\infty)$, and suppose that
$X$ is lattice-isomorphic to a  sublattice of $Y(I)$,
where $I=[0,1]$ or $[0,\infty).$ Then there is a lattice
embedding $V:X\to Y(\Omega\times [0,\infty))$ of the form
$Vx(\omega,t)=a(\omega,t)x(\omega)$ where $a$ is a
nonnegative Borel function on $\Omega\times [0,\infty)$ of
the form:  $$ a(\omega,t)=\sum_{k\in\bold
Z}2^{m(k,\omega)}e_k(t)$$ where $m:\bold
Z\times\Omega\to\bold Z\cup\{-\infty\}$ is a Borel map such
that $k\to m(k,\omega)$ is decreasing for each $\omega.$
Furthermore if $I=[0,1]$ then $a$ is supported on a set of
finite measure.  \endproclaim

If $Y$ is an r.i. space on $I$=$[0,1]$ or $[0,\infty)$ and
$f\in Y_+\setminus\{0\}$ then we define $Y_f$ to be the r.i.
space on $I$ defined by $y\in Y_f$ if and only if $y\otimes
f\in Y(I\times I)$ where $y\otimes f(s,t)=y(s)f(t).$ The
norm on $Y_f$ is given by $\|y\|_{Y_f}=\|y\otimes f\|_Y.$
Notice that since $f$ dominates a function of the form
$\alpha\chi_E$ where $\alpha>0$ and $\lambda(E)>0$ there
exists a constant $C$ depending on $f$ so that $\|y\|_Y \le
C\|y\|_{Y_f}.$

\proclaim{Proposition 3.3}Suppose $Y$ is an order-continuous
quasi-Banach r.i. space on $[0,\infty)$ and that $X$ is an
order-continuous quasi-Banach r.i. space on $[0,1]$.  Let
$U:X\to Y$ be a lattice homomorphism and let
$U\chi_{[0,1]}=f\neq 0$.  Then:\newline (1) There exists $C$
so that if $x\in X$ then $\|x\|_{Y_f}\le C\|x\|_X.$\newline
(2) If $U$ is a lattice embedding then $X=Y_f[0,1].$
\endproclaim

\demo{Remark}If $U$ is a lattice embedding of $X$ into
$Y[0,1]$ then the above proposition gives $X=Y_f[0,1]$ where
$f\in Y[0,1].$\enddemo

\demo{Proof}We use Lemma 3.1 to construct the lattice
homomorphism $V:X\to Y([0,1]\times [0,\infty)).$ Notice that
if $g\in Y[0,\infty)$ has the same distribution as
$V\chi_{[0,1]}$ then $Y_f[0,1]=Y_g[0,1]$ with equivalent
norms.

Let $u$ be any nonnegative simple function on $[0,1]$ of the
form $u=\sum_{j=1}^n\alpha_j\chi_{B_j}$ where
$\{B_1,\ldots,B_n\}$ is a Borel partition of $[0,1].$ For
any $N$ let $$a_N(s,t)=\sum_{|k|\le N}\sum_{|m(k,s)|\le
N}2^{m(k,s)}e_k(t)$$ and let $b_N=a-a_N.$ We can partition
$g=g_N+h_N$ where $g_N$ has the same distribution as $a_N$
and $h_N$ has the same distribution as $b_N.$

Now $a_N =\sum_{|k|\le N}\sum_{|l|\le
N}2^l\chi_{A_{kl}}(s)e_k(t)$ where $(A_{kl})_{k,l}$ are
Borel subsets of $[0,1].$ We can therefore use Liapunoff's
theorem to find Borel sets $B'_1,\ldots,B'_n$ so that
$\lambda(B'_j)=\lambda(B_j)$ for all $j$ and $\lambda(B'_j
\cap A_{kl})=\lambda(B_j)\lambda(A_{kl})$ whenever $1\le
j\le n$ and $-N\le k,l\le N.$ Let
$u'=\sum_{j=1}^n\alpha_j\chi_{B'_j}$.  Then $a_N(s,t)u'(t)$
has the same distribution as $u\otimes g_N$.  Hence $$
\|u\otimes g_N\|_Y \le \|Vu'\|_Y \le \|u\otimes g_N\|_Y +
\|u\|_{\infty}\|\chi_{[0,1]}\otimes h_N\|_Y.$$

For case (1) we let $N\to \infty$ and deduce that
$\|u\|_{Y_g}\le \|U\|\|u\|_X.$

For case (2) we observe that, since $Y$ is order-continuous,
$$\lim_{N\to\infty}\|\chi_{[0,1]}\otimes h_N\|_Y=0.$$ Since
$V$ is an embedding there exists $c>0$ so that we have a
lower-estimate $\|Vu'\|_Y\ge c\|u\|_X$.  Hence $\|u\|_X \le
c^{-1}\|u\|_{Y_g}.$

If $X$ lattice embeds into $Y[0,1]$ then $a$ has support of
measure at most one and hence so has $f$ so that we can
assume that $f\in Y[0,1]$.  \bull\enddemo

\proclaim{Corollary 3.4} Suppose $Y$ is an order-continuous
quasi-Banach r.i. space on $[0,\infty)$ and that $X$ is an
order-continuous quasi-Banach r.i. space on $[0,1]$.  Let
$U:X\to Y$ be a lattice homomorphism.  If $U\neq 0$ then
there is a constant $C$ so that $\|x\|_Y\le C\|x\|_X$ for
$x\in X[0,1].$\endproclaim

\demo{Proof}This follows from (1) of the preceding
proposition combined with the remarks before
it.\bull\enddemo

\demo{Remark}This Corollary is well-known (see Abramovich
\cite{1} and remarks in the introduction to \cite{27}).
\enddemo

For our final result of this section, we will need the
following factorization theorem, which is essentially due to
Krivine \cite{31} (\cite{34}); we will, however, prove the
form of the theorem required here.

\proclaim{Proposition 3.5}Suppose $0< p<\infty$.  Suppose
$Y$ is an $p$-concave quasi-K\"othe function space on
$(\Omega,\mu)$ and suppose that either (a)
$P:L_p(\Delta,\lambda)\to Y$ is a lattice homomorphism or
(b) $p\ge 1$ and $P:L_p(\Delta,\lambda)\to Y$ is a positive
operator.  Then there is a Borel function $w\in L_0(\mu)$
with $w>0$ a.e. so that $$ \|f\|_Y \le \|fw\|_p$$ for $f\in
L_0(\mu)$ and $$ \|w(Pf)\|_p\le \|P\|\|f\|_p$$ for $f\in
L_p(\Delta).$ \endproclaim

\demo{Proof}We can suppose $P\neq 0.$ We require the
following property of $P$ which is valid in cases (a) or
(b):  if $f_1,\ldots,f_n\ge 0$ in $L_p$ then
$P((\sum_{i=1}^nf_i^p)^{1/p}) \ge
(\sum_{i=1}^n(Pf_i)^p)^{1/p}$ (see \cite{34} p. 55).  Let
$u$ be any strictly positive function in $Y$.  Now consider
the subsets $E$ and $F$ of $L_{\infty}$ defined by
$E=\{f:f\ge 0,\ \|uf^{1/p}\|_Y> \|P\|\}$ and $F=\{f:\exists
\ 0\le x\in L_p, \|x\|_p\le 1 ,\ u^pf\le (Px)^p\}$.

It is clear that $E$ is convex.  We argue that $\text{co }F$
does not meet $E.$ Indeed suppose $f_1,\ldots,f_n\in F$ and
$c_1,\ldots,c_n\ge 0$ with $\sum_{j=1}^nc_j=1.$ Suppose
$u^pf_j\le (Px_j)^p$ where $x_j\ge 0$ and $\|x_j\|_p\le 1.$
Then $u^p(\sum_{j=1}^nc_jf_j)\le \sum_{j=1}^nc_j(Px_j)^p\le
(Py)^p$ where $y=(\sum_{j=1}^nc_jx_j^p)^{1/p}$ so that
$\|y\|_p\le 1$ (see \cite{34} Proposition 1.d.9).  Since $F$
includes the negative cone it has non-empty interior.  Now,
by the Hahn-Banach theorem, there exists $\Phi\in
L_{\infty}^*$ so that $\Phi(f-g)> 0$ if $f\in E$ and $g\in
F.$ Clearly $\Phi\ge 0,$ and $\Phi(f)>0$ if $f\ge 0$ and
$f\neq 0;$ hence since $P$ is not zero we have $\inf_{g\in
E}\Phi(g)>0.$ By normalizing we can suppose $\inf_{g\in
E}\Phi(g)=1.$ Let us write $\Phi(f)=\int
f\phi\,d\mu+\Phi_0(f)$ where $\phi\in L_1(\mu),$ and
$\Phi_0$ is singular with respect to $\mu.$ If $f\in E$ we
may find $0\le f_n\uparrow f$ a.e. so that
$\Phi(f_n)\uparrow \int f\phi\,d\mu.$ However by order
continuity $f_n\in E$ for large enough $n$ and so $\int
f\phi\,d\mu\ge 1$ for $f\in E.$

Now it is clear that if $y\in Y_+$ with $\|y\|_Y =1.$ Then
for $\epsilon>0$ we have that
$(\|P\|+\epsilon)^py^{p}u^{-p}\in E$ and so
$\|yu^{-1}\phi^{1/p}\|_p\ge \|P\|^{-1}.$ Thus if $y\in Y$
then $\|y\|_Y\le \|yw\|_p$ where $w=\|P\|\phi^{1/p}u^{-1}.$
If $f\in L_p(\Delta,\lambda)$ with $\|f\|_p=1$ then
$(P(|f|))^pu^{-p}\in F$ and so $$\int (P(|f|))^p\phi
u^{-p}d\mu \le 1,$$ so that $\|wP(|f|)\|_p \le \|P\|$ which
implies the theorem.  \bull\enddemo

\proclaim{Theorem 3.6}Suppose $0<p<\infty$ and $Y$ is a
$p$-concave quasi-Banach r.i. space on $[0,1]$ or
$[0,\infty)$.  Suppose $L_p$ is lattice-isomorphic to a
sublattice of $Y.$ Then $Y[0,1]=L_p[0,1].$ \endproclaim
\demo{Proof} It suffices to consider the case when
$Y=Y[0,\infty).$ By Proposition 3.3 there exists $f\in Y$ so
that $Y_f[0,1]=L_p[0,1].$ Thus there is a lattice embedding
$V:L_p\to Y([0,1]\times[0,\infty))$ of the form $x\to
x\otimes f.$ We assume $\|x\|_p\le \|x\otimes f\|_Y\le
C\|x\|_p.$ Applying Proposition 3.5, there is a nonnegative
weight function $w$ on $[0,1]\times [0,\infty)$ so that
$\|y\|_Y \le \|yw\|_p$ for $y\in Y$ and $\|x\|_p\le
\|x\otimes f\|_Y\le \|w(x\otimes f)\|_p \le C\|x\|_p$ for
$x\in L_p.$

Now let $v(t)=(\int_0^1w(s,t)^pds)^{1/p}.$ It follows from a
symmetrization argument that if $y\in Y$ then $$ \|y\|_Y \le
\left(\int_0^1\int_0^{\infty}v(t)^p|y(s,t)|^p
ds\,dt\right)^{1/p},$$ and that $$
\int_0^{\infty}f(t)^pv(t)^pdt \le C^p.$$

Let $u$ be the increasing rearrangement of $v$ so that
$u(t)=\inf_{\lambda(E)=t}\sup_{s\in E}v(s).$ Then if as
usual $y^*$ is the decreasing rearrangement of $|y|$, the
first equation yields that if $y\in Y[0,\infty)$ then $$
\|y\|_Y \le (\int_0^{\infty}y^*(t)^pu(t)^p dt)^{1/p}.$$

In particular for $0<s<1,$ $$ s \le \|D_sf^*\|_Y^p \le
\int_0^{\infty}f^*(t/s)^pu(t)^pdt.$$ This in turn implies
that $$ \int_0^{\infty}f^*(t)^pu(st)^pdt \ge 1.$$ Now $\int_0^{\infty}
f^*(t)^pu(t)^pdt\le C^p.$ Letting $s\to 0$ we obtain from
the Dominated Convergence Theorem that $\lim_{t\to
0}u(t)=c>0$ and $\int_0^{\infty}f(t)^pdt \le C^pc^{-p}.$

Pick $0<\tau<\infty$ so that $\|f^*\chi_{[\tau,\infty)}\|_Y
\le 1/2.$ It follows from $p$-concavity that
$$\|D_s(f^*\chi_{[\tau,\infty)})\|_Y \le s^{1/p}/2.$$ On the
other hand $\|D_sf^*\|_Y=\|\chi_{[0,s]}\otimes f\|_Y \ge
s^{1/p}.$ Hence $\|D_s(f^*\chi_{[0,\tau]})\|_Y \ge
s^{1/p}/2.$ From this and $p$-concavity we also obtain
easily that
$$(\frac1{s\tau}\int_0^{s\tau}f^*(t/s)^pdt)^{1/p}\|\chi_{[0,s\tau)}\|_Y
\ge \frac12 s^{1/p}.$$ Hence $\|\chi_{[0,t]}\|_Y \ge
c_1t^{1/p}$ when $0\le t\le 1$ for a suitable constant
$c_1.$ This in turn implies, by $p$-concavity, that if $y\in
Y[0,1]$ then $\|y\|_Y\ge c_1\|y\|_p$ and this is enough to
show that $Y[0,1]=L_p[0,1].$ \bull\enddemo

\vskip10pt \heading{4.  Cone-embeddings}\endheading
\vskip10pt

Let $X$ and $Y$ be quasi-Banach lattices.  We will say that
a positive operator $L:X\to Y$ is a {\it cone-embedding} if
$L$ satisfies a lower bound for positive elements, i.e.
there exists $\delta>0$ so that $\|Lx\|_Y\ge \delta\|x\|_X$
for $x\ge 0.$ We will say that $L$ is a {\it strong
cone-embedding} if it additionally satisfies the condition
that for some $C>0$ and every $x_1,\ldots,x_n\ge 0$ we have
$\|\max_{1\le k\le n}x_k\|_X \le C\|\max_{1\le k\le
n}Lx_k\|_Y.$ This is trivially equivalent to requiring the
same inequality for $x_1,\ldots,x_n$ mutually disjoint.

Our first results demonstrate conditions under which every
cone-embedding is a strong cone-embedding.

\proclaim{Lemma 4.1}Suppose $s,\delta>0,$ and
$1<p,q<\infty$.  Then there is a constant
$C=C(s,p,q,\delta)$ so that if $X$ is a $p$-convex K\"othe
function space, $Y$ is an $s$-convex, $q$-concave
quasi-K\"othe function space (where each constant of
convexity and concavity is one) and if $L:X\to Y$ is a
cone-embedding satisfying $ \delta\|x\|_X \le \|Lx\|_Y \le
\|x\|_X$ for $x\ge 0$ then if $x_1,\ldots,x_n\ge 0$ are
disjoint, $$ \|\sum_{j=1}^nx_j\|_X \le C\|\max_{1\le j\le
n}Lx_j\|_Y.$$ \endproclaim

\demo{Proof}We pick $m=m(p,\delta)$ so that
$2^{m(1-1/p)}\delta>2.$

First notice that if $x_1,\ldots,x_n$ are disjoint, $$\left(
\Ave_{\epsilon_{ij}=\pm 1} \sum_{j=1}^n
\prod_{i=1}^m(1+\epsilon_{ij})^px_j^p\right)^{1/p}=2^{m(1-1/p)
}\sum_{j=1}^n x_j.$$ Thus by $p$-convexity $$
2^{m(1-1/p)}\|\sum_{j=1}^nx_j\|_X \le \left( \Ave_{\epsilon_{ij}=\pm 1}
\|\sum_{j=1}^n\prod_{i=1}^m(1+\epsilon_{ij})x_j\|_X^p\right)^{1/p}.$$

Now it follows that $$ \align
2^{m(1-1/p)}\|\sum_{j=1}^nx_j\|_X &\le\delta^{-1}\left(
\Ave_{\epsilon_{ij}=\pm1}
\|\sum_{j=1}^n\prod_{i=1}^m(1+\epsilon_{ij})Lx_j\|_Y^p\right)^{1/p}\\
&\le
\delta^{-1}\sum_{I\subset[m]}\left(\Ave_{\epsilon_{ij}=\pm1}
\|\sum_{j=1}^n\prod_{i\in
I}\epsilon_{ij}Lx_j\|_Y^p\right)^{1/p}\\ &\le
\delta^{-1}\left(\|\sum_{j=1}^nLx_j\|_Y
+C_1(2^m-1)\|(\sum_{j=1}^n|Lx_j|^2)^{1/2}\|_Y \right),
\endalign $$ where $C_1=C_1(q,s),$ using Theorem 1.d.6 of
\cite{34}.

Reorganizing we have, since $2^{m(1-1/p)}\delta-1>1,$ and
$\|L\|\le 1,$ $$ \align \delta^{-1}\|\sum_{j=1}^nx_j\|_X
&\le C_12^m\|(\sum_{j=1}^n|Lx_j|^2)^{1/2}\|_Y\\ &\le C_12^m
\|\sum_{j=1}^n Lx_j\|_Y^{1/2}\|\max_{1\le j\le
n}Lx_j\|_Y^{1/2} \endalign $$ and this in turn implies,
since $Y$ is $s$-convex for some $s>0,$ $$
\|\sum_{j=1}^nx_j\|_X \le C_1^22^m\delta^2\|\max_{1\le j\le
n}Lx_j\|_Y.\bull$$\enddemo

Let us give a simple application.

\proclaim{Theorem 4.2}Suppose $Y$ is an $r$-convex Banach
lattice where $r>2$ which is $q$-concave for some
$q<\infty.$ Suppose that $X$ is a $p$-convex Banach lattice,
where $p>2,$ which is isomorphic to a subspace of $Y.$ Then
$X$ is $r$-convex.\endproclaim

\demo{Remarks}This result is well-known for $1< r\le 2$  (cf.
\cite{34}, p. 51).  The hypothesis on
$X$
is equivalent to the statement that $\ell_2^n$ is not lattice finitely
representable in $X$ (note that $X$ must be of type 2, and
apply Lemma 2.4 of \cite{21}).  In \cite {21} there are two
results closely related to Theorem 4.2.  Theorem 2.3 of
\cite{21} is the analogous result for upper $r$-estimates in
place of $r$-convexity, while Theorem 2.6 (or Proposition
2.e.10 of \cite{34}) implies the above theorem for the
special case when $X$ is an r.i. space on $[0,1]$.  In this
latter case one can replace the hypothesis that $X$ is
strictly 2-convex by the weaker hypothesis that $X\neq
L_2[0,1].$ \enddemo

\demo{Proof}It suffices to consider the case when the
$r$-convexity, $q$-concavity constants of $Y$ are both one
and the $p$-convexity constant of $X$ is one.  We may also
suppose that $X$ and $Y$ are K\"othe function spaces.  We
will suppose that there is a bounded linear operator $S:X\to
Y$ with $\delta\|x\|_X\le \|Sx\|_Y\le \|x\|_X.$ It will also
suffice to prove the result when $X$ is finite-dimensional,
i.e.  $\Omega=\{1,2,\ldots,n\}$ and thus has a
1-unconditional basis $(e_k)_{k=1}^n$ consisting of atoms,
provided we establish a uniform bound on the $r$-convexity
constant $M^r(X)$ in terms of $(p,q,r,\delta).$

To this end we define a map $L:X_{1/2}\to Y_{1/2}$ by
$Le_k=|Se_k|^2.$ It follows from Krivine's theorem that if
$x=\sum_{k=1}^n\xi_ke_k\ge 0$ then $$ \align
\|Lx\|_{Y_{1/2}}&=
\|(\sum_{k=1}^n\xi_k|Se_k|^2)^{1/2}\|_Y^2\\ &\le
K_G^2\|(\sum_{k=1}^n\xi_ke_k)^{1/2}\|_X^2\\ &\le K_G^2
\|x\|_{X_{1/2}}.  \endalign $$ Also since $Y$ is
$q$-concave, there exists $C_0=C_0(q)$ so that $$ \align
\|Lx\|_{Y_{1/2}} &\ge
C_0^{-2}(\Ave_{\epsilon_k=\pm1}\|\sum_{k=1}^n\epsilon_k
\xi_k^{1/2}Se_k\|_Y)^{2}\\ &\ge
C_0^{-2}\delta^{-2}\|x\|_{X_{1/2}}.  \endalign $$ Now by
Lemma 4.1 applied to $K_G^{-2}L$, using the fact that
$X_{1/2}$ is $p/2$-convex and $Y_{1/2}$ is $r/2-$convex and
$q/2-$concave we obtain the existence of
$C_1=C_1(p,q,r,\delta)$ so that for
$x=\sum_{k=1}^n\xi_ke_k\ge 0$, $$ \|x\|_{X_{1/2}}\le C_1
\|\max_{1\le k\le n} \xi_k Le_k\|_{Y_{1/2}}$$ which in turn
implies that if $x\in X,$ with $x=\sum\xi_ke_k,$ $$ \|x\|_X
\le C_2\|\max_{1\le k\le n}|\xi_k||Se_k|\|_Y$$ where
$C_2^2=C_1.$ Now suppose $x_1,\ldots,x_m\in X$ with
$x_j=\sum_{k=1}^n\xi_{jk}e_k.$ Then $$ \align
\|(\sum_{j=1}^m|x_j|^r)^{1/r}\|_X &\le C_2\|\max_{1\le k\le
n}(\sum_{j=1}^m|\xi_{jk}|^r)^{1/r}|Se_k|\|_Y \\ &\le C_2
\|(\sum_{k=1}^n\sum_{j=1}^m|\xi_{jk}|^r|Se_k|^r)^{1/r}\|_Y
\\ &= C_2
\|(\sum_{j=1}^m(\sum_{k=1}^n|\xi_{jk}|^r|Se_k|^r))^{1/r}\|_Y\\
&\le C_2 (\sum_{j=1}^m
\|(\sum_{k=1}^n|\xi_{jk}|^r|Se_k|^r)^{1/r}\|_Y^r )^{1/r}\\
&\le C_2 (\sum_{j=1}^m
\|(\sum_{k=1}^n|\xi_{jk}|^2|Se_k|^2)^{1/2}\|_Y^r )^{1/r}\\
&\le
K_GC_2(\sum_{j=1}^m\|(\sum_{k=1}^n|\xi_{jk}|^2|e_k|^2)^{1/2}\|_X^r
)^{1/r}\\ &\le K_GC_2 (\sum_{j=1}^m\|x_j\|_X^r)^{1/r}.
\endalign $$ This completes the proof.\bull\enddemo

We now give a second criterion for a cone-embedding to be a
strong cone-embedding.

\proclaim{Lemma 4.3}Suppose $0<q,s<\infty$ and that $X$ is
an $s$-convex quasi-Banach r.i. space on $[0,1]$ or
$[0,\infty)$ with $p_X>1.$ Suppose $Y$ is an $s$-convex
$q$-concave quasi-K\"othe function space and $L:X\to Y$ is a
cone-embedding.  Then there is a constant $C$ so that if
$x_1,\ldots,x_n\ge 0$ are disjoint, $$ \|\sum_{j=1}^nx_j\|_X
\le C\|\max_{1\le j\le n}Lx_j\|_Y.$$ \endproclaim

\demo{Proof}We suppose that $\|L\|\le 1$ and that $\delta>0$
is such that if $x\ge 0$ then $\delta\|x\|_X\le \|Lx\|_Y\le
\|x\|_X.$ We may also suppose that for some $p>1$ and some
constant $C_0$ we have $\|D_t\|_X\le C_0t^{1/p}$ for $t\ge
1.$

We select first an integer $m$ so that $2^{m(p-1)} \ge
2^{p+1}C_0^p\delta^{-p}.$ Let $\theta=2^{-m}.$

Now suppose $x_1,\ldots,x_n\ge 0$ are given; it will suffice
to consider the case when each $x_i$ is a countably simple
function (i.e. takes only a countable set of values) and
$\|\sum_{i=1}^nx_i\|_X=1.$ Suppose $N$ is an integer with
$N>4(2^mn).$ Then for each $1\le i\le n$ we can write
$x_i=\sum_{j=1}^Nx_{ij}$ as a disjoint sum where
$x_{ij}^*=D_{(1/N)}x_i^*$.

Let $\epsilon_{ijk}=\pm 1$ be a choice of signs for $1\le
i\le n,\ 1\le j\le N$ and $1\le k\le m$ and denote by
$\epsilon$ the array $(\epsilon_{ijk}).$ We define $$
u(\epsilon)=\sum_{i=1}^n\sum_{j=1}^N
\prod_{k=1}^m(1+\epsilon_{ijk})x_{ij}.$$

Let $\xi_i(\epsilon)$ be the number of $j$ such that
$\epsilon_{ijk}=1$ for $1\le k\le m.$ As functions on the
natural finite probability space of all choices of signs
$\epsilon$, the functions $\xi_i$ for $1\le i\le n$ are
independent and identically distributed with binomial
distributions corresponding to a sample size $N$ and
probability for an individual trial of $\theta.$ They each
have mean $\alpha=N\theta$ and variance
$N\theta(1-\theta)<\alpha.$ Notice that by choice of $N$ we
have $\alpha>4n$

We thus have $$ \Ave_{\epsilon}\max_{1\le i\le
n}|\xi_i-\alpha|^2 \le \sum_{i=1}^n\Ave_{\epsilon}
|\xi_i-\alpha|^2 < n\alpha.$$ Let
$\zeta(\epsilon)=\min_{1\le i\le n}\xi_i(\epsilon).$ Then $$
\Ave_{\epsilon}|\alpha-\zeta|^2 < n\alpha$$ and so $$
\Ave_{\epsilon} \zeta > \alpha-(n\alpha)^{1/2}
>\frac12\alpha.$$

We next turn to estimating $\Ave \|u(\epsilon)\|_X^p.$ In
fact we have that for each $\epsilon,$ $$ \|\sum_{i=1}^n
x_i\|_X \le 2^{-m}\|D_{N/\zeta}u(\epsilon)\|_X.$$ Thus we
have an estimate that $$ 1\le C_0\theta
N^{1/p}\zeta(\epsilon)^{-1/p}\|u(\epsilon)\|_X.$$
Reorganizing and averaging gives $$ \Ave_{\epsilon}
\|u(\epsilon)\|_X^p \ge
C_0^{-p}\theta^{-p}N^{-1}\Ave_{\epsilon} \zeta \ge \frac12
C_0^{-p}\theta^{1-p}.$$ The original choice of $m$ now gives
the estimate $$ \Ave_{\epsilon}\|u(\epsilon)\|_X^p \ge
2^p\delta^{-p}$$ which implies $$\left (\Ave_{\epsilon}
\|L(u(\epsilon))\|_Y^p\right)^{1/p} \ge 2.$$

We now proceed as in Lemma 4.1, expanding out and concluding
that for some constant $C_1$ depending only on $Y,$ $$
\left(\Ave_{\epsilon}\|L(u(\epsilon))\|_Y^p \right)^{1/p}
\le \|\sum_{i=1}^n\sum_{j=1}^NLx_{ij}\|_Y + 2^mC_1
\|(\sum_{i=1}^n\sum_{j=1}^N |Lx_{ij}|^2)^{1/2}\|_Y.$$ Since
$\|\sum_{i=1}^n\sum_{j=1}^NLx_{ij}\|_Y \le 1$ we can
conclude that $$ \|(\sum_{i=1}^n\sum_{j=1}^N
|Lx_{ij}|^2)^{1/2}\|_Y \ge C_1^{-1}\theta.$$ Again this
implies that $$ \|\max_{1\le i\le n}Lx_i\|_Y \ge
\|\max_{i,j}Lx_{ij}\|_Y \ge C_1^{-2}\theta^{-2}.$$ The
result now follows.\bull\enddemo

\proclaim{Proposition 4.4}Let $X$ be an order-continuous
K\"othe function space on $(\Delta,\lambda),$ which contains
$L_{\infty}.$ Let $Y$ be a quasi-K\"othe function space on
$(\Omega,\mu)$ which is $s$-convex for some $s>0$ and
$q$-concave for some $q<\infty.$ Suppose $L:X\to Y$ is a
strong cone-embedding.  Then, for $n\ge 1,$ there exist
Borel maps $a_n:\Omega\to [0,\infty)$ and
$\sigma_n:\Omega\to \Delta$ so that $a_1\ge a_2\ge\cdots\ge
0$ and $\sigma_m(\omega)\neq\sigma_n(\omega)$ if $m\neq n,$
and for some $C>0$ we have for any $x\in X$ with $x\ge 0,$
$$ C^{-1}\|x\|_X \le \|\max_na_nx\circ\sigma_n\|_Y \le
\|\sum_{n=1}^{\infty}a_nx\circ\sigma_n\|_Y \le C\|x\|_X.$$
\endproclaim

\demo{Proof}We use the random measure representation of
positive operators (see \cite{25},\cite{41},\cite{42}).
There exists a Borel map $\omega\to\nu_{\omega}$ from
$\Omega$ to $\Cal M(\Delta),$ endowed with the weak$^*$
topology, so that for any $x\in X$ we have $$ Lx(\omega)
=\int x(t)d\nu_{\omega}\ \ \ \mu-\text{a.e.}$$ Further we
can write
$$\nu_{\omega}=\sum_{n=1}^{\infty}a_n(\omega)\delta_{\sigma_n(\omega)}+
\nu'_{\omega}\ \ \ \mu-\text{a.e}$$ where $a_n:\Omega\to
[0,\infty)$ and $\sigma_n:\Omega\to\Delta$ are Borel maps
satisfying the assumptions above, and $\nu'_{\omega}$ is a
continuous measure.

Since $L$ is a strong cone-embedding there exists a constant $C$
so that $\|L\|\le C$ and whenever $x_1,\ldots,x_n$ are
disjoint and positive in $X$ then $$ \|\sum_{j=1}^nx_j\|_X
\le C\|\max_{1\le j\le n} Lx_j\|_Y.$$

Now suppose $x\ge 0$.  Then for each $m,$ $$ \|x\|_X \le
C\|\max_{E\in\Cal A_m}L(x\chi_E)\|_Y.  $$ For the definition
of $\Cal A_m$ see Section 2. Now $\max_{E\in\Cal
A_m}L(x\chi_E)$ is monotone decreasing to
$\max_na_nx\circ\sigma_n$ so that, by the order-continuity
of $Y,$ $$ C^{-1}\|x\|_X\le \|\max_na_nx\circ\sigma_n\|_Y
\le \|\sum_{n=1}^{\infty}a_nx\circ\sigma_n\|_Y\le
\|Lx\|_Y\le C\|x\|_X.\bull$$ \enddemo

\demo{Remark}Of course there is no special significance in
modelling $X$ on $(\Delta,\lambda)$ here; we clearly have
the same result for any Polish measure space $(K,\nu)$.
Note also that in the above argument the pointwise maximum
$\max_n a_nx\circ\sigma_n$ exists $\mu-$a.e. for $x\in X.$
\enddemo

\proclaim{Proposition 4.5}Suppose $Y$ is an order-continuous
quasi-Banach r.i. space on $[0,\infty)$ with property (d).
Suppose that either $X$ is an order-continuous atomic
quasi-Banach lattice or that $X$ is an order-continuous
quasi-K\"othe function space on $(\Delta,\lambda)$ and that
$L:X\to Y$ is a strong cone-embedding.  Then $X$ is
lattice-isomorphic to a sublattice of $Y.$\endproclaim

\demo{Proof}Let $C$ be a constant greater than the property
(d) constant of $Y$ and the constant in the definition of
the strong cone-embedding.  Let us prove this first for the
case when $X$ is atomic.  Then we regard $X$ as a sequence
space (a quasi-K\"othe space modelled on $\bold N$).  Let
$(e_n)_{n\in\bold N}$ be the basis vectors and let
$u_n=Le_n.$ We define a map $V:X\to L_0(\bold N\times
[0,\infty))$ by $Ve_n=v_n$ where $v_n(k,t)=0$ if $k\neq n$
and $v_n(n,t)=u_n(t).$ If $a_1,\ldots,a_n\ge 0$ then it is
easy to see that $$ (\sum_{k=1}^na_kv_k)^{**}\le
(\sum_{k=1}^na_ku_k)^{**}$$ and so by property (d) we have
that $V$ is bounded and $\|V\|\le C\|L\|.$ However since $L$
is a strong cone-embedding $$ \|\sum_{k=1}^na_ke_k\|_X \le
C\|\max_{1\le k\le n}a_ku_k\|_Y \le
C\|\sum_{k=1}^na_kv_k\|_Y$$ so that $V$ is an isomorphism
onto its range.

The nonatomic case is similar.  We can suppose that $X$ is a
quasi-K\"othe function space on $(\Delta,\lambda)$
containing $L_{\infty}$ and that $L$ is of the form $$ Lx
=\sum_{n=1}^{\infty}a_nx\circ\sigma_n$$ where for some
constant $C_1$ we have $$C_1^{-1}\|x\|_X \le
\|\max_na_nx\circ\sigma_n\|_Y \le
\|\sum_{n=1}^{\infty}a_nx\circ\sigma_n\|_Y \le C_1\|x\|_X.$$
Define $V:X \to L_0(\bold N\times [0,\infty))$ by the
formula $ Vx(n,t) = a_n(t)x(\sigma_n(t)).$ Then if $x\ge 0$
we have $(Vx)^{**}\le (Lx)^{**}$ so that $V$ is bounded,
while $$C_1^{-1}\|x\|_X \le \|\max_na_nx\circ\sigma_n\|_Y
\le \|Vx\|_Y,$$ so that $V$ is also an isomorphism.
\bull\enddemo

\proclaim{Proposition 4.6} Suppose $Y$ is an
order-continuous quasi-Banach r.i. space on $[0,1]$ with
property (d).  Suppose for some $p>1,$ $Y_{1/p}$ has
property (d).  Suppose $X$ is an order-continuous
quasi-K\"othe function space on $[0,1],$ and that $L:X\to Y$
is a strong cone-embedding.  Then there is a Borel subset
$E$ of $[0,1]$ with $\lambda(E)>0$ so that $X(E)$ is
lattice-isomorphic to a sublattice of $Y.$\endproclaim

\demo{Proof}We again may suppose that $X$ is a quasi-K\"othe
function space containing $L_{\infty}$.  Note first that we
must have $Y_{1/p}\subset L_1$ and hence $Y\subset L_p.$ We
may extend $Y$ to be a quasi-Banach r.i. space on
$[0,\infty)$ in several different ways.  Precisely we define
$W$ to be the space of $f\in L_0[0,\infty)$ so that
$f^*\chi_{[0,1]}\in Y$ and $f\in L_1[0,\infty)$ with the
associated quasi-norm
$\|f\|_{W}=\max(\|f^*\chi_{[0,1]}\|_Y,\|f\|_1).$ We define
$Z$ to be the space of $f\in L_0[0,\infty)$ so that
$f^*\chi_{[0,1]}\in Y$ and $f\in L_p[0,\infty)$ with the
associated quasi-norm
$\|f\|_{Z}=\max(\|f^*\chi_{[0,1]}\|_Y,\|f\|_p).$ Then both
$W$ and $Z$ have property (d).  Note that both $W[0,1]$ and
$Z[0,1]$ coincide with $Y$ and hence $L$ may be regarded as
mapping into either $W$ or $Z$.  Note also that $W\subset Z$
with continuous inclusion.

Appealing to the preceding Proposition, we can find a
lattice embedding $U:X\to W[0,\infty)$ in such a way that
for some constant $C$ we have $ C^{-1}\|x\|_X \le
\|Ux\|_{Z}$ and $\|Ux\|_{W}\le C\|x\|_X$ for $x\ge 0.$

It now follows by Lemma 3.1 and Proposition 3.2 that we can
find a a nonnegative Borel function $a$ on $[0,1]\times
[0,\infty)$ with $a(t,s))$ decreasing in $s$ for each fixed
$t$ so that the map $Vx(t,s)=a(t,s)x(t)$ defines a lattice
embedding of $X$ into $Z_1([0,1]\times [0,\infty))$ and such
that for some $C_1$ we have $ C_1^{-1}\|x\|_X \le
\|Vx\|_{Z}$ and $\|Vx\|_{W}\le C_1\|x\|_X$ for $x\ge 0.$

Notice in particular that $$
\int_0^1\int_0^{\infty}a(t,s)^rds\,dt <\infty$$ for $r=1$
and $r=p.$ We therefore can find constants $0<c<M<\infty$ so
that there is a Borel subset $E$ of $[0,1]$ of positive
measure such that if $t\in E$ then $
\int_0^{\infty}a(t,s)dt\le M$ and
$\int_0^{\infty}a(t,s)^p\ge c^p.$ If $t\in E$ then
$a(t,s)\le Ms^{-1}$ and so we also have
$\int_{s_0}^{\infty}a(t,s)^pds \le M^ps_0^{1-p}.$

Recall that $Y_{1/p}$ has property (d) and therefore
$Z_{1/p}$ also has property (d) and is an interpolation
space between $L_1$ and $L_{\infty}$ with some constant
$\gamma^p\ge 1.$ Pick $u>1$ so that $cu^{1-1/p} >4\gamma M.$
We then modify $V$ to form $V_0:X\to Z$ by setting
$V_0x=Vx\chi_{E\times [0,u]}.$ We will show that $V_0$ is a
lattice embedding of $X(E)$ into $Z([0,1]\times
[0,\infty)).$

Let $P$ be the positive operator defined on $L_1([0,1]\times
[0,\infty) )$ and $L_{\infty}([0,1]\times [0,\infty))$ by $$
Pg(t,s)
=\left(\frac1u\int_0^ug(t,v)dv\right)a(t,s)^p\chi_E(t)\chi_{[u,\infty)}(s).$$
It is easy to calculate that $\|Pg\|_{\infty}\le
M^pu^{-p}\|g\|_{\infty}.$ Similarly $\|Pg\|_1\le
M^pu^{1-p}.$ It follows that $\|P\|_{Z_{1/p}} \le
M^pu^{1-p}\gamma^p.$

It follows that if $f\in Z$ then $$ \|(P(|f|^p))^{1/p}\|_Z
\le \gamma Mu^{\frac1p-1}\|f\|_Z\le \frac{12}c\|f\|_Z.$$

Suppose in particular $x\in X(E)$ and $x\ge 0.$ Let $f=Vx.$
Then $$
P(|f|^p)(t,s)=x(t)\left(\frac1u\int_0^ua(t,v)^pdv\right)^{1/p}
a(t,s)\chi_E(t)\chi_{[u,\infty)}(s).$$ However for $t\in E$
$$ \int_0^ua(t,v)^pdv \ge c^p- M^pu^{1-p}\ge \frac12c^p\ge
(c/2)^p.$$ Hence $$ \|Vx-V_0x\|_Z \le 2\gamma M
u^{\frac1p-1} c^{-1}\|Vx\|_Z\le \frac12\|Vx\|_Z.$$ It
follows that $V_0$ maps $X(E)$ isomorphically into
$Z([0,1]\times [0,u] )$ which is lattice isomorphic to $Y.$
\bull\enddemo

\proclaim{Corollary 4.7} Suppose $Y$ is an order-continuous
quasi-Banach r.i. space on $[0,1]$ with property (d).
Suppose for some $p>1,$ $Y_{1/p}$ has property (d).  Suppose
$X$ is an order-continuous quasi-Banach r.i. space on
$[0,1],$ and that  $L:X\to Y$ is a strong
cone-embedding.  Then there exists $f\in Y_+\setminus\{0\}$
such that $X=Y_f.$ \endproclaim

\demo{Proof}This follows from Proposition 3.3.\bull\enddemo

\vskip10pt \heading{5.  Cone-embeddings of r.i.
spaces}\endheading \vskip10pt

\proclaim{Proposition 5.1}Suppose $0<s<q<\infty$ and that
$X$ is an $s$-convex, $q$-concave quasi-K\"othe function
space on $(\Omega,\mu).$ Suppose $m\ge q$ is a natural
number.  Then there is a constant $C=C(X)$ so that if
$x_1,\ldots,x_n\in X_+$ and $b_1,\ldots,b_n\ge 0$ then $$
(\Ave_{\pi\in\Pi_n}\|\sum_{i=1}^nb_{\pi(i)}x_i\|^m)^{1/m}
\le C\max\left((\Ave_{\pi\in \Pi_n}\|\max_{1\le i\le
n}b_{\pi(i)}x_i\|^m)^{1/m},\
\frac1n(\sum_{i=1}^nb_i)\|\sum_{i=1}^nx_i\|\right).$$
\endproclaim

\demo{Proof}This is a somewhat disguised form of the
so-called Classification Formula (Theorem 2.1 of \cite{21}
or Theorem 2.e.5 of \cite{34}).  It can be derived from this
formula; we indicate the direct proof.  We assume that $X$
has $q$-concavity constant one.  Then $Z=L_m(\Pi_n:X)$ has
$m$-concavity constant one where $\Pi_n$ is given its
natural probability measure.  Now there is a constant $C_0$
depending only on $m$ so that if $f_1,\ldots,f_n\in Z_+$ $$
(\sum_{i=1}^nf_i)^m\le C_0\left(\sum_{|A|=m}(\prod_{i\in
A}f_i) +
(\sum_{i=1}^nf_i^2)(\sum_{i=1}^nf_i)^{m-2})\right).$$ Hence
for $C_1=C_1(s,m)$ $$ \|\sum_{i=1}^nf_i\|_Z \le
C_1\max(\|S_1\|,\|S_2\|)$$ where $S_1=
(\sum_{|A|=m}(\prod_{i\in A}f_i))^{1/m}$ and $S_2=(\max
f_i)^{1/m}(\sum f_i)^{1-1/m}.$ Now as $Z$ is $s$-convex we
can estimate:  $$ \align \|S_2\| &\le (\|\max
f_i\|)^{1/m}(\|\sum f_i\|)^{1-1/m} \\ &\le \frac{1}{m}\|\max
f_i\| + (1-\frac1m)\|\sum f_i\| \endalign $$ It then follows
that $$ \|\sum f_i\| \le \max (m\|S_1\|,\|\max f_i\|).$$

Now let $f_i=b_i\xi_i$ where $\xi_i(\pi)=x_{\pi(i)}$.  Then
if $n\ge 2m,$ we use $m$-concavity:  $$ \align \|S_1\| &=
(\Ave_{\pi}\|(\sum_{|A|=m}\prod_{i\in
A}b_ix_{\pi(i)})^{1/m}\|^m)^{1/m}\\ &\le
\|(\Ave_{\pi}\sum_{|A|=m}\prod_{i\in
A}b_ix_{\pi(i)})^{1/m}\|\\ &\le
(\frac{(n-m)!}{n!})^{1/m}(\sum_{i=1}^nb_i)\|\sum_{i=1}^nx_i\|\\
&\le \frac2n(\sum_{i=1}^nb_i)\|\sum_{i=1}^nx_i\|.  \endalign
$$ The Proposition now follows easily.\bull\enddemo

\proclaim{Proposition 5.2}Let $X,Y$ be order-continuous
quasi-Banach r.i. spaces on $[0,1]$.  Suppose that $p_Y>1$,
$Y$ is $q$-concave for some $q<\infty$ and that there is a
cone-embedding $L:X\to Y$.  Then either $X=L_1[0,1]$ or $X$
is lattice-isomorphic to a sublattice of $Y$ and so $X=Y_f$
for some $f\in Y_+.$\endproclaim

\demo{Proof}For ease of notation we regard $X$ as modelled
on $(\Delta,\lambda).$

Let us first note that the proof is trivial if we assume
$p_X>1.$ Indeed in this case $L$ is a strong cone-embedding
(Lemma 4.3 ) and $Y_{1/r}$ has property (d) as long as
$1<r<p_Y.$ So Corollary 4.7 applies.  We therefore need only
to prove that if $X\neq L_1$ then $p_X>1.$

Assume then $X\neq L_1.$ Note first that $Y\subset L_r$ if
$1<r<p_Y.$

We can assume that, for some $\delta>0$ and every $x\in X,$
$\delta\|x\|_X\le \|Lx\|_Y \le \|x\|_X$ for $x\ge 0.$ Let us
consider the random measure representation of $L$ i.e.
$$Lx(s) =\int x\,d\mu_s$$ where $s\to\mu_s$ is a
weak$^*$-Borel map from $[0,1]$ to $\Cal M(\Delta).$ We can
as usual write
$$\mu_s=\sum_{n=1}^{\infty}a_n(s)\delta_{\sigma_n(s)}
+\nu_s$$ where $a_n:[0,1]\to [0,\infty)$ and
$\sigma_n:[0,1]\to\Delta$ are Borel maps and
$\sigma_m(s)\neq\sigma_n(s)$ if $m\neq n,$ and $\nu_s$ is
for each $s$ nonatomic.

Since $Y$ has nontrivial concavity there is a constant $C_0$
and an integer $m$ so that if $y_1,\ldots,y_n\in Y_+$ and
$b_1,\ldots,b_n\ge 0,$ $$
(\Ave_{\pi\in\Pi_n}\|\sum_{i=1}^nb_{\pi(i)}y_i\|_Y)^{1/m}
\le C_0\max\left((\Ave_{\pi\in \Pi_n}\|\max_{1\le i\le
n}b_{\pi(i)}y_i\|_Y^m)^{1/m},\
\frac1n(\sum_{i=1}^nb_i)\|\sum_{i=1}^ny_i\|_Y\right).\tag
*$$

Let us introduce the functional on $X$ defined by
$$\Gamma(x)= \sup\{ \|\max_na_nu\circ\sigma_n\|_Y:\
u^*=x^*\}.$$

Consider a nonnegative simple function $x\in
CS_{n_0}(\Delta).$ For each $n\ge n_0$ we can write
$x=\sum_{E\in\Cal A_n}\xi_E\chi_E.$ For each permutation
$\pi$ of $\Cal A_n$ let $x_{\pi}=\sum_{E\in\Cal
A_n}\xi_{\pi(E)}\chi_E.$ Let $y_E=L\chi_E\in Y.$

We also define for each $n,$ and each $s\in [0,1]$
$\tau_n(s)$ to be the least integer $\tau$ so that
$(\sigma_i(s))_{i=1}^{\tau}$ belong to distinct members of
$\Cal A_n.$ Note that $\lim_{n\to\infty}\tau_n(s)=\infty$
for all $s.$

Note that $$ \max_{E\in\Cal A_n}\xi_{\pi(E)}y_E(s) \le
\max_{1\le k\le \tau_n}a_kx_{\pi}\circ\sigma_k
+\|x\|_{\infty}\left(\sum_{k>\tau_n}a_k(s)+\max_{E\in\Cal
A_n}\nu_s(E)\right),$$ so that $$ \|\max_{E\in\Cal
A_n}\xi_{\pi(E)}y_E\|_Y \le \Gamma(x)
+\eta_n\|x\|_{\infty}$$ where $\lim_{n\to\infty}\eta_n=0.$
Now appealing to (*) gives that $$ \delta\|x\|_X \le
C_0\max(\Gamma(x),\|x\|_1\|L\chi_{\Delta}\|_Y)
+C_0\eta_n\|x\|_{\infty}$$ which gives us $$ \delta\|x\|_X
\le C_0\max(\Gamma(x),\|x\|_1\|L\chi_{\Delta}\|_Y).  \tag
**$$

First suppose $a_1$ vanishes a.e. so that $\Gamma(x)=0$ for
all $x\ge 0.$ Then $$\|x\|_X \le
C_0\delta^{-1}\|L\chi_{\Delta}\|_Y\|x\|_1$$ for all $x\in X$
so that $L_1\subset X.$ Since $X\neq L_1$ we must have that
$X^*=\{0\}$ and Theorem 4.4 of \cite{25} shows that $L$ must
vanish (we remark that in the preparatory Lemma 4.3 of
\cite{25} the hypothesis $X^*=\{0\}$ has been omitted in the
statement).  This is impossible so we must have that $a_1>0$
on a set of positive measure.  Hence if we set
$Sx=a_1x\circ\sigma_1$ then $S$ is a nontrivial lattice
homomorphism of $X$ into $Y$ and Corollary 3.3 will yield
that $X\subset Y$.  Hence $X\subset L_r$ where $1<r<p_Y.$

We next show that in fact $\|x\|_X \le C_1\Gamma(x).$ If
not, there is a sequence $x_n$ with $\|x_n\|_X=1,$ $x_n\ge
0$ and $\Gamma(x_n)\to 0.$ But, if this happens we must have
$x_n\to 0$ in measure and $\|x_n\|_r$ bounded.  Hence
$\lim_{n\to\infty}\|x_n\|_1=0$ and $(**)$ yields that
$\lim_{n\to\infty}\|x_n\|_X=0.$ This contradiction
establishes the claim.

Fix any simple $f\in X[0,1]$ with $\|f\|_X=1$.  Then there
exists $0\le x \in X(\Delta)$ with $x^*=f^*$ so that $$
\|\max_na_nx\circ\sigma_n\|_X \ge C_1^{-1}.$$ Let
$x=\sum_{j=1}^M\xi_j\chi_{H_j}$ where $H_1,\ldots,H_M$ are
disjoint Borel sets.  For each $s$ let $k(s)$ be the first
index such that $a_k(s)x(\sigma_k(s))=\max_{1\le
n<\infty}a_n(s)x(\sigma_n(s)).$ Then let $b'(s)=a_{k(s)}(s)$
and $\rho(s)=\sigma_{k(s)}(s).$ The operator $V:X\to Y$
given by $Vz=b'z\circ\rho$ is then a lattice homomorphism
form $X$ into $Y$ with $\|V\|\le 1.$ For $n\in\bold Z$ let
$F_n=(b')^{-1}(2^n,2^{n+1}]$ and let $b=\sum_{n\in\bold
Z}2^n\chi_{F_n}$ so that $\frac12b'\le b\le b'.$ For each
$n$ the measures $B\to \lambda(\rho^{-1}B\cap F_n)$ are
absolutely continuous.  Then for any $N,$ we can use
Liapunoff's theorem to find sets $H_j^{\alpha}\subset H_j$
with $\lambda(H_j^{\alpha})=\alpha^{-1} \lambda(H_j)$ for
$1\le j\le M$ and $\lambda(\rho^{-1}H_j^{\alpha}\cap
F_n)=\alpha^{-1}\lambda(\rho^{-1}H_j\cap F_n)$ for $|n|\le
N.$ Let $G_N=\cup_{|n|\le N}F_n.$ Then $$
\|V(\sum_{j=1}^M\xi_j\chi_{H_j}^{\alpha})\|_Y \ge
\frac12\|D_{1/\alpha}(\chi_{G_N}Vx)\|_Y\ge
\frac12\|D_{\alpha}\|_Y^{-1}\|\chi_{G_N}Vx\|_Y.$$ Letting
$N\to \infty$ we have $$ \|Vx\|_Y \le2 \|D_{\alpha}\|_Y
\|D_{\alpha^{-1}}x\|_X.$$ Hence $$ \|f\|_X \le
2C_1\|D_{\alpha}\|_Y \|D_{\alpha}^{-1}f\|_X.$$ As this
inequality holds for all simple $f\ge 0$ we obtain $$
\|D_{\alpha}\|_X \le 2C_1 \|D_{\alpha}\|_Y$$ for $\alpha>1$
so that $p_X\ge p_Y>1$.  As observed in the introductory
remarks, this is sufficient to prove the
theorem.\bull\enddemo

The following Proposition is trivially false in the case
when $p=1$ since the map $x\to
(\int_0^1x(s)\,ds)\chi_{[0,1]}$ is a cone embedding of
$L_1[0,1]$ into $L_p[0,1]$ when $p<1.$

\proclaim{Proposition 5.3}Suppose $1<p<\infty.$ Suppose $Y$
is a $p$-concave quasi-Banach r.i. space on $[0,1]$ or
$[0,\infty)$ and that there is a cone-embedding of
$L_p(\Delta,\lambda)$ into $Y$.  Then $Y[0,1]=L_p[0,1].$
\endproclaim

\demo{Proof} We assume that $Y$ is $s$-normed.  Let
$(\Omega,\mu)$ represent either $[0,1]$ or $[0,\infty)$ with
associated Lebesgue measure.  We apply Lemma 4.1 and
Proposition 4.4.  There exists a constant $C$ and Borel maps
$a_n:\Omega\to [0,\infty)$ and $\sigma_n:\Omega\to \Delta$
so that $\sigma_m(\omega)\neq \sigma_n(\omega)$ if $m\neq n$
and so that for $0\le x\in L_p,$ $$ C^{-1}\|x\|_p \le
\|\max_n a_nx\circ\sigma_n\|_Y \le
\|\sum_{n=1}^{\infty}a_nx\circ\sigma_n\|_Y \le C\|x\|_p.$$
Let $Lx=\sum a_nx\circ\sigma_n;$ then $L:L_p\to Y$ is a
positive operator.  We can also apply Proposition 3.5:
there is a weight function $w>0$ on $\Omega$ so that
$\|y\|_Y \le \|wy\|_p$ for $y\in Y$ and $\|w(Lx)\|_p\le
C\|x\|_p$ for $x\in L_p.$

At this point we define measures $\nu_n$ on $\Delta$ by
$\nu_n(B)=\int_{\sigma_n^{-1}B}w^pa_n^pd\mu.$ It is easy to
see that each $\nu_n$ is a finite Borel measure absolutely
continuous with respect to $\lambda.$ Hence we can find
derivatives $v_n=d\nu_n/d\lambda.$ Now if $0\le x\in
L_p(\Delta)$ then $$
\int_{\Omega}w^p\sum_{n=1}^{\infty}a_n^p(x\circ\sigma_n)^p
d\mu=\int_{\Delta} x^p(\sum_{n=1}^{\infty}v_n)d\lambda$$ and
so it follows that $$\sum_{n=1}^{\infty} v_n(t) \le C^p$$
almost everywhere.  By an application of Egoroff's theorem
we can find a Borel set $E\subset \Delta$ of positive
measure and $N$ so that $$ \sum_{n=N+1}^{\infty} v_n(t) \le
(2^{1/s}C)^{-p}$$ for $t\in E.$

Now, observe that if $0\le x\in L_p(E)$ then $$ \align
\|\max_{n\ge N+1}a_nx\circ\sigma_n\|_Y &\le
\|w(\sum_{n=N+1}^{\infty}a_n^p(x\circ\sigma_n)^p)^{1/p}\|_p\\
&\le (\int_E x^p(\sum_{n=N+1}^{\infty}v_n)d\lambda)^{1/p}\\
&\le \frac1{2^{1/s}C}\|x\|_p.  \endalign $$ Hence $$
\|\max_{1\le n\le N}a_nx\circ\sigma_n\|_Y \ge
(\|Lx\|_Y^s-\frac1{2C^s}\|x\|_p^s)^{1/s}\ge
\frac1{2^{1/s}C}\|x\|_p.$$ This implies that $L_p$ is
isomorphic to a sublattice of $Y^N$ and hence to a
sublattice of $Y.$ Finally we can apply Theorem 3.6 to
deduce that $Y[0,1]=L_p[0,1].$ \bull\enddemo

\vskip10pt \heading{6.  The main construction}
\endheading\vskip10pt

\proclaim{Lemma 6.1}Let $X$ be a $q$-concave K\"othe
function space on some Polish measure space $(\Omega,\mu)$,
where $q<\infty.$ Then there is a constant $C$ depending
only on $X$ so that if $f_1,\ldots,f_n\in X,$ and
$h=(\sum_{i=1}^nf_i^2)^{1/2},$ then for any $M>1$ we have:
$$ (\Ave_{\epsilon_i=\pm
1}\|g_{\epsilon}\chi_{H_{\epsilon}}\|^q)^{1/q} \le
CM^{-1}\|h\| $$ where
$g_{\epsilon}=\sum_{i=1}^n\epsilon_if_i$ and
$H_{\epsilon}=\{|g_{\epsilon}|\le
M^{-1}h\}\cup\{|g_{\epsilon}|\ge Mh\}.$\endproclaim

\demo{Proof}Note first that $$
\Ave_{\epsilon_i=\pm1}\|g_{\epsilon}\chi_{(|g_{\epsilon}|\le
M^{-1}h)}\|^q \le M^{-q}\|h\|^q.$$ On the other hand, if
$C_0$ is the $q$-concavity constant of $X,$
$$(\Ave_{\epsilon_i=\pm1}
\|g_{\epsilon}\chi_{(|g_{\epsilon}|\ge Mh)}\|^q)^{1/q} \le
C_0^{-1}\|\phi\|$$ where $$ \phi(s)=
(\int_{g_{\epsilon}(s)\ge
Mh(s)}|\sum_{i=1}^m\epsilon_if_i(s)|^qd\epsilon)^{1/q}.$$ We
can estimate (assuming $h(s)>0$) $$ \phi(s)^q \le
M^{-q}h(s)^{-q}\int|\sum_{i=1}^m\epsilon_if_i(s)|^{2q}d\epsilon
\le C_1^qM^{-q}h(s)^q $$ where $C_1$ is a constant
determined by the constant in Khintchine's inequality for
$2q.$ Combining we have $$ (\Ave_{\epsilon_i=\pm1}
\|g_{\epsilon}\chi_{(|g_{\epsilon}|\ge Mh)}\|^q)^{1/q} \le
C_0^{-1}C_1M^{-1}\|h\|.$$ The result now follows.\bull
\enddemo

We now introduce some notation.  If $[a,b]$ is a closed
interval with $1<a$ we write $\Gamma(a,b)$ for the
collection of measurable functions $f$ on $[0,\infty)$ which
satisfy that almost everywhere, either $f(s)=0$ or $a\le
|f(s)|\le b$ or $b^{-1}\le |f(s)|\le a^{-1}.$ Let
$[a_n,b_n]_{n=0}^{\infty}$ be a sequence of intervals with
$a_0=1.$ If $(\eta_n)_{n=0}^{\infty}$ is a sequence with
$0<\eta_n<1$ then $[a_n,b_n]$ is {\it $(\eta_n)$-separated}
if $b_n\le \eta_na_{n+1}$ for all $n.$

\proclaim{Lemma 6.2}Let $X$ be an r.i. space on
$[0,\infty).$ Suppose $0<\delta<1,$ and that $\sigma>0$ is
such that $2^5\sigma <\delta.$ Suppose
$(\eta_n)_{n=1}^{\infty} $ is any sequence satisfying $\sum
\eta_n<\sigma$ and that $[a_n,b_n]_{n=0}^{\infty}$ are
$(\eta_n)-$separated then for any $f_0,\ldots, f_N\in X$
such that \newline (1) $\delta\le \|f_j\| \le 1$ for $0\le
j\le N$\newline (2) $f_j\in \Gamma(a_j,b_j)$ \newline we
have that $(f_j)_{j=0}^{N}$ is 2-equivalent to a disjointly
supported sequence $(g_j)_{j=0}^N$ with $g_j\in
\Gamma(a_j,b_j)$ for $0\le j\le N.$ \endproclaim

\demo{Proof} Let $E_k=\{s:|f_k(s)|=\max_{0\le j\le
N}|f_j(s)|,\ |f_k(s)|>|f_j(s)|\ \text{ if }j<k\}.$ Let
$g_k=f_k\chi_{E_k}.$ Then (with appropriate modifications if
$k=0$ or $k=N$) $$ \align \|f_k-g_k\| &\le
\sum_{j<k}\|f_k\chi_{(|f_k|\le |f_j|)}\|
+\sum_{j>k}\|f_k\chi_{(|f_k|<|f_j|)}\| \\ &\le
\sum_{j<k}a_k^{-1}b_j\|f_j\| +\sum_{j>k}b_ka_j^{-1}\|f_j\|
\\ &\le \sum_{j<k}\prod_{i=j}^{k-1}\eta_i +
\sum_{j>k}\prod_{i=k}^{j-1}\eta_i\\ &\le
(\eta_{k-1}+\eta_k)\prod_{i=1}^{\infty}(1+\eta_i)\\ &\le
e^{\sigma}(\eta_{k-1}+\eta_k)\le 4(\eta_k+\eta_{k-1}).
\endalign $$ Hence $\|g_k\|\ge \delta-8\sigma\ge
\frac{\delta}2.$ We also have $\sum\|f_k-g_k\| \le 8\sigma
\le \frac{\delta}4.$ Since $(g_k)$ is a disjoint sequence it
follows from standard perturbation theory that $(f_k)$ is
2-equivalent to $(g_k).$ \bull\enddemo

\proclaim{Lemma 6.3} Let $X$ be an r.i. space on
$[0,\infty)$ or $[0,1].$ Suppose $0<\delta<\frac12$ and
$[a_n,b_n]_{n=0}^{\infty}$ are
$(2^{-(n+6)}\delta)-$separated.  Then for any positive
disjoint $f_0,f_1,\ldots,f_N\in X$ such that $\delta\le
\|f_j\|\le 1$ and $f_j\in \Gamma(a_j,b_j)$ for $0\le j\le N$
we have that $(f_j)_{j=0}^N$ is $6$-equivalent to a
disjointly supported sequence in $E_X.$ \endproclaim

\demo{Proof}We suppose at first that $X$ is an r.i. space on
$[0,\infty).$ For $0\le j\le N$ we choose $m_j\in\bold Z$ so
that $$ \frac{\delta}{2^{j+4}b_j}<\|\chi_{[0,2^{m_j}]}\| \le
2\frac{\delta}{2^{j+4}b_j}.$$ Similarly we choose
$n_j\in\bold Z$ for $1\le j\le N$ so that $$
\frac{a_j\delta}{2^{j+4}}<\|\chi_{[0,2^{n_j}]}\|\le
2\frac{a_j\delta}{2^{j+4}}.$$ It is clear that $m_N\le
m_{N-1}\le\cdots\le m_0\le n_1\le\cdots\le n_N.$

Let $u_j=\lambda(a_j\le f_j\le b_j)$ and
$v_j=\lambda(b_j^{-1}\le f_j\le a_j^{-1}).$ Then if $1\le
j\le N,$ $$ \align \|\chi_{[0,2^{m_j}+u_j]}\| &\le
\|\chi_{[0,2^{m_j}]}\|+\|\chi_{[0,u_j]}\|\\ &\le
2^{-(j+3)}b_j^{-1}\delta + a_j^{-1}\\ &\le 2a_j^{-1}\\ &\le
2^{-(j+4)}\delta b_{j-1}^{-1}\\ &<\|\chi_{[0,2^{m_{j-1}}]}\|
\endalign $$ so that $2^{m_j}+u_j\le 2^{m_{j-1}}.$
Similarly, if $1\le j\le N-1,$ $$ \|\chi_{[0,2^{n_j}+v_j]}\|
\le 2^{-(j+3)}a_j\delta + b_j<2b_j$$ so that $$
\|\chi_{[0,2^{n_j}+v_j]}\| \le 2^{-(j+5)}\delta a_{j+1}
<\|\chi_{[0,2^{n_{j+1}}]}\|$$ and $2^{n_j}+v_j <
2^{n_{j+1}}.$ Finally $$ \|\chi_{[0,2^{m_0}+u_0+v_0]}\| \le
2^{-3}\delta b_0^{-1} +b_0<2b_0\le \frac12\delta a_1.$$
Hence $2^{m_0}+u_0+v_0<2^{n_1}.$

It now follows that we can rearrange $f_0,\ldots,f_N$ in the
following manner.  We can suppose that $f_0$ is supported
and decreasing on $[2^{m_0},2^{n_1}).$ Let
$f'_0=f_0\chi_{[2^{m_0},1]}$ and $f''_0=f_0-f'_0.$ For $1\le
j\le N,$ we let $f'_j=f_j\chi_{(a_j\le f_j\le b_j)}$ and
$f''_j =f_j\chi_{(b_j^{-1}\le f_j\le a_j^{-1})}.$ We can
then suppose that for $1\le j\le N,$ $f'_j$ is supported and
decreasing on $[2^{m_j},2^{m_{j-1}})$ and $f''_j$ is
supported and decreasing on $[2^{n_j},2^{n_{j+1}})$ where we
adopt the convention $n_{N+1}=\infty.$

Now if $e_k=\chi_{[2^k,2^{k+1}]}$ let $$x'_0
=\sum_{k=m_0}^{-1}f'_0(2^{k+1})e_k$$ and $$ x''_0
=\sum_{k=0}^{n_1-1}f''_0(2^{k+1})e_k.$$ and for $1\le j\le
N,$ let $$ x'_j =\sum_{k=m_j}^{m_{j-1}-1}f'_j(2^{k+1})e_k$$
and $$x''_j=\sum_{k=n_j}^{n_{j+1}-1}f''_j(2^{k+1})e_k.$$ We
set $x_j=x'_j+x''_j$ for $0\le j\le N.$

Then $0\le x_j\le f_j$ for $0\le j\le N$.  However if
$D_2g(t)=g(t/2)$ we have $f_j \le D_2x_j +z_j$ where $ z_j=
b_je_{m_j} + a_j^{-1}e_{n_j}$ for $1\le j\le N$ and $z_0\le
b_0e_{m_0}.$ Thus $\|z_j\|_X \le 2^{-(j+2)}\delta$ for $0\le
j\le N.$ Hence if $\alpha_0,\ldots,\alpha_N\ge 0,$ $$ \align
\|\sum_{j=0}^N\alpha_jf_j\|&\le 2\|\sum_{j=0}^N
\alpha_jx_j\| + \sum_{j=0}^N2^{-(j+2)}\delta\alpha_j \\ &\le
2\|\sum_{j=0}^N\alpha_jx_j\| +
\frac{\delta}2\max|\alpha_j|\\ &\le
2\|\sum_{j=0}^N\alpha_jx_j\|+\frac12\|\sum_{j=0}^N\alpha_jf_j\|
\endalign $$ so that $(f_j)_{j=0}^N$ is 4-equivalent to a
disjoint sequence in $E_X.$ This completes the proof when
$X$ is modelled on $[0,\infty).$

For the case $X=X[0,1]$, we may regard $X$ as being defined
on $[0,\infty)$ and proceed as before, but with each $f_j$
having support of measure at most one.  In this case, we
have $x''_0=0$ while for $1\le j\le N,$ we have $x''_j\le
f''_j\le a_j^{-1} \le 2^{-6j}\delta.$ Hence if $\alpha_j\ge
0$ for $0\le j\le N, $ $$ \|\sum_{j=0}^N\alpha_jx''_j\| \le
2^{-5}\delta\max|\alpha_j|\le
\frac{1}{32}\|\sum_{j=0}^N\alpha_jf_j\|.$$ Hence since
$\frac12+\frac1{32}<\frac23,$ $$ \align
\|\sum_{j=0}^N\alpha_jf_j\|&\le
2\|\sum_{j=0}^N\alpha_jx_j\|+\frac12\|\sum_{j=0}^N\alpha_jf_j\|
\\ &\le
2\|\sum_{j=0}^N\alpha_jx'_j\|+\frac23\|\sum_{j=0}^N\alpha_jf_j\|
\endalign $$ and $(f_j)_{j=0}^N$ is 6-equivalent to
$(x'_j)_{j=0}^N$ which is a sequence in
$E_{X[0,1]}$.\bull\enddemo

We now consider a situation which will remain fixed for
Lemmas 6.4-6.6.  We suppose now that $X$ is a good K\"othe
function space on $(\Delta,\lambda)$ which is $q$-concave
with constant one where $q<\infty.$ We further suppose that
$Y$ is an r.i. space on $[0,\infty)$ which is also
$q$-concave with constant one.  We will assume that $X$ is
isomorphic to a subspace of $Y.$ Let us therefore suppose
that $T:X\to Y$ is a bounded linear operator satisfying
$\delta \|x\|_X\le \|Tx\|_Y\le \|x\|_X$ where $\delta>0.$

For convenience we recall the notation introduced in Section
2. For $\epsilon_k=\pm1,$ we denote by
$\Delta(\epsilon_1,\ldots,\epsilon_n)$ the clopen subset of
$\Delta$ of all $(d_j)_{j=1}^{\infty}$ such that
$d_j=\epsilon_j$ for $1\le j\le n.$ For each $n$ let $\Cal
A_n$ denote the collection of
$\Delta(\epsilon_1,\ldots,\epsilon_n)$.  Let $\Cal C_n$ be
the algebra generated by the atoms $\Cal A_n.$ We let $CS_n$
denote the linear span of $\{\chi_E:E\in\Cal A_n\}.$ We also
define the Haar functions
$h_E=\chi_{\Delta(\epsilon_1,\ldots,\epsilon_n,+1)}-
\chi_{\Delta(\epsilon_1,\ldots,\epsilon_n,-1)}$ for
$E=\Delta(\epsilon_1,\ldots,\epsilon_n).$ Let $CS$ be the
union of the spaces $CS_n.$

We define $Q_n:CS_n\to L_0[0,\infty)$ to be the linear map
such that $Q_n(\chi_E)=|Th_E|^2$ where $E\in\Cal A_n.$

\proclaim{Lemma 6.4}If $x\in CS_n$ then
$\|(Q_nx^2)^{1/2}\|_Y \le K_G \|x\|_X$, where $K_G$ is the
Grothendieck constant.\endproclaim

\demo{Proof} If $x=\sum_{E\in \Cal A_n}\alpha_E\chi_E$ then,
by Krivine's theorem \cite{31}, $$\|(\sum
|\alpha_E|^2|Th_E|^2)^{1/2}\|_Y \le K_G \| \left ( \sum
|\alpha_E|^2|h_E|^2\right)^{1/2}\|_X=K_G\|x\|_X.\bull$$\enddemo

For any measurable function $f\in L_0[0,\infty)$ and $a\ge
1$ we define $\tau_af=f\chi_{(a^{-1}\le |f|\le a)}.$ We then
define for $x\in CS_+,$ $$\Psi(x)=
\sup_a\liminf_{m\to\infty}\|\tau_a((Q_m(x^2))^{1/2})\|_Y.$$

\proclaim{Lemma 6.5}There exists a constant $C=C(X,Y)$ so
that if $\eta>0$ and $b\ge 1,$ then whenever $x\ge 0,$ $x\in
CS_n$ with $\Psi(x)<\eta$ then there exists a clopen set $D$
independent of $\Cal C_n$ such that $\lambda(D)=\frac12,$
$\max(\|x\chi_D\|_X,\|x-x\chi_D\|_X)\le (3/4)^{1/q}\|x\|_X$
and $$ \|\tau_b(T(x-2\chi_D x))\| \le C\eta.  $$
\endproclaim

\demo{Proof}Let $\eta=\Psi(x).$ We first pick $w\in
L_0(\Delta)_+$ so that $\|x\|_X=\|xw^{-1}\|_q$ and
$\|\xi\|_X \le \|\xi w^{-1}\|_q$ for all $\xi\in X.$ We
write $x=\sum_{E\in\Cal A_n}\alpha_E \chi_E.$

Suppose $m\ge n$.  For a choice of signs $\epsilon_F=\pm1$
we write $$x_{\epsilon}= \sum_{E\in\Cal A_n}
\alpha_E\sum_{F\in\Cal A_m\atop F\subset E} \epsilon_F
h_F.$$ We also let $y_{\epsilon}=Tx_{\epsilon}\in Y.$

Let $x_{\epsilon,+}=\max(x_{\epsilon},0)$ and
$x_{\epsilon,-}=\max(-x_{\epsilon},0).$ We first estimate $$
\|x_{\epsilon,+}\|_X^q\le \frac12\int
|x|^qw^{-q}\sum_{F\in\Cal A_m}(\chi_F+\epsilon_F
h_F)d\lambda.$$ This gives $$
\|x_{\epsilon,+}\|_X^q-\frac12\|x\|_X^q \le \left|\sum_{F\in
\Cal A_m} \epsilon_F\int_F|x|^qw^{-q}h_Fd\lambda\right|.$$
Switching signs we get a similar estimate for
$\|x_{\epsilon,-}\|_X^q$ and hence $$
\Ave_{\epsilon_i=\pm1}\max(\|x_{\epsilon,+}\|_X^q,\|x_{\epsilon,-}\|_X^q)\le
\frac12\|x\|_X^q +\left(\sum_{F\in \Cal A_m}\left(\int_F
|x|^qw^{-q}d\lambda)\right)^2\right)^{1/2}$$ by Khintchine's
inequality.

The second term here can be estimated by $$ \max_{F\in\Cal
A_m}\left(\int_F|x|^qw^{-q}d\lambda\right)^{1/2}\|x\|_X^{q/2}.$$
It follows that for large enough $m$ we have $$
\Ave_{\epsilon_i=\pm1}\max(\|x_{\epsilon,+}\|_X^q,\|x_{\epsilon,-}\|_X^q)
\le \frac58\|x\|_X^q.$$ For such $m$ we have $$\Pr(\max
(\|x_{\epsilon,+}\|_X^q,\|x_{\epsilon,-}\|_X^q) \le
\frac34\|x\|^q)\ge \frac16.$$

We will now choose $m$ subject to this restriction and such
that $$ \|\tau_a(Q_mx^2)^{1/2}\|_Y \le \eta$$ where
$a=b\|x\|_X/\eta.$ Let $G=\{a^{-1}\le (Q_mx^2)^{1/2}\le
a\}.$ Then since $Y$ has cotype $q$, for a suitable constant
$C_0=C_0(Y),$ $$ \align \left(\Ave_{\epsilon_i=\pm1}
\|y_{\epsilon}\chi_G\|_Y^q\right)^{1/q} &\le C_0\|\chi_G(
\sum_{E\in\Cal A_n}\sum_{F\in \Cal A_m\atop F\subset
E}|\alpha_E|^2|Th_E|^2)^{1/2}\|_Y \\ &= C_0 \|\chi_G
(Q_mx^2)^{1/2}\|_Y \\ &\le C_0\eta.  \endalign $$ On the
other hand, if $H$ is the complement of $G$ and
$B_{\epsilon}=\{b^{-1}\le |y_{\epsilon}|\le b\}$ then
$B_{\epsilon}\cap H \subset \{|y_{\epsilon}|\le
\eta\|x\|_X^{-1} (Q_mx^2)^{1/2}\}\cup
\{|y_{\epsilon}|\ge\|x\|_X \eta^{-1}(Q_mx^2)^{1/2}\}.$ It
thus follows from Lemmas 6.1 and 6.4 that $$
\left(\Ave_{\epsilon_i=\pm1}\|\tau_by_{\epsilon}\chi_H\|_Y^q\right)^{1/q}
\le C_1\eta\|x\|_X^{-1}\|(Q_mx^2)^{1/2}\|_Y \le K_GC_1
\eta.$$

Hence $$
\left(\Ave_{\epsilon_i=\pm1}\|\tau_by_{\epsilon}\|_Y^q
\right)^{1/q} \le C_2\eta$$ where $C_2$ depends only on
$X,Y.$

Finally it follows there must exist a choice of $\epsilon_F$
so that $\max(\|x_{\epsilon,+}\|^q,\|x_{\epsilon,-}\|^q) \le
\frac34 \|x\|_X^q$ and $\|\tau_by_{\epsilon}\|_Y\le
6C_2\eta.$ We conclude by writing $\sum\epsilon_Fh_F =
2\chi_D-\chi_{\Delta}$ and then $D$ satisfies our
hypotheses.\bull\enddemo

\proclaim{Lemma 6.6}Suppose $\inf\{\Psi(x):\|x\|_X=1,\ x\in
CS_+\}=0.$ Then there is a nonatomic Banach lattice $Z$
which is lattice-finitely representable in $X$ so that $Z$
has an unconditional basis which is lattice-finitely
representable in $E_Y.$ \endproclaim

\demo{Proof}Suppose $N$ is a natural number.  Let
$\gamma=(\frac34)^{1/q}.$ Let $C$ be the constant determined
in the previous lemma.  We will select $\eta>0$ so that $$
\eta <
\min(\frac12\gamma^N\delta,\frac{(1-\gamma)^N\delta}{10^22^{N+1}(C+1)}).$$

We pick $x\in CS_+$ so that $\|x\|_X=1$ and $\Psi(x)<\eta.$
Suppose $x\in CS_n.$ We construct by induction a sequence of
clopen sets $(F_k)_{k=1}^{2^{N+1}-1}$, sequences
$(a_k,b_k)_{k=0}^{2^{N}-1}$, and functions $y_k \in Y$ for
$0\le k\le 2^{N+1}-1$ so that:\newline (1) $a_0=1$ and
$F_1=\Delta.$ \newline (2) Each $F_k$ is independent of
$\Cal C_n.$\newline (3) $F_k=F_{2k}\cup F_{2k+1}$ and
$\lambda(F_{2k})=\lambda(F_{2k+1})=\frac12\lambda(F_k)$ for
$1\le k\le 2^N-1.$\newline (4) For $1\le k\le 2^N-1$ we have
$(1-\gamma)\|x\chi_{F_k}\|_X \le
\|x\chi_{F_{2k}}\|_X,\|x\chi_{F_{2k+1}}\|_X \le
\gamma\|x\chi_{F_k}\|_X.$ (5) $a_k\le b_k$ ($1\le k\le
2^{N+1}-1$) and $b_k \le (1-\gamma)^N\delta
2^{-(k+7)}a_{k+1}$ for $1\le k \le 2^{N+1}-2.$\newline (6)
If $h_0=x$ and then $h_k=x(2\chi_{F_{2k}}-\chi_{F_k})$ for
$1\le k\le 2^N-1$ then $\|Th_k-y_k\|<(C+1)\eta.$\newline (7)
$y_k\in \Gamma(a_k,b_k).$

We start the induction as stated with
$a_0=1,F_1=\Delta,h_0=x.$ We then select $b_0$ large enough
so that $\|Th_0-\tau_{b_0}Th_0\|_Y<\eta$ and set
$y_0=\tau_{b_0}Th_0.$

Now suppose $1\le k\le 2^N-1$ and that $(a_j)_{j=0}^{k-1}$,
$(b_j)_{j=0}^{k-1}$, $(y_j)_{j=0}^{k-1}$ and
$(F_j)_{j=1}^{2k-1}$ have been determined.  We first pick
$a_k$ so that $b_{k-1}\le (1-\gamma)^N\delta 2^{-(k+6)}a_k$
so that (5) holds.  Now $\Psi(x\chi_{F_k})\le\Psi(x)<\eta.$
Hence we are able to apply Lemma 6.5 to find a clopen set
$D$ independent of the algebra generated by the sets $\Cal
C_n$ and $\{F_1,\ldots,F_{k-1}\}$ so that
$\lambda(D)=\frac12$, $\max(\|x\chi_{F_k\cap
D}\|_X,\|x\chi_{F_k\setminus D}\|_X) \le
\gamma\|x\chi_{F_k}\|_X,$ and
$$\|\tau_{a_k}(T(x\chi_{F_k}-2x\chi_{F_k\cap D}))\|_Y \le
C\eta$$ where $C$ is the constant of the previous lemma.

We now let $F_{2k}=F_k\cap D$ and $F_{2k+1}=F_k\setminus D.$
Conditions (2) and (3) are immediately satisfied.  Condition
(4) follows from the triangle law.  If we define $h_k$ by
(6) we have $\|\tau_{a_k}Th_k\|_Y \le C\eta.$ Therefore we
can pick $b_k>a_k$ so large that if $G$ is the set where
$b_k^{-1}\le |Th_k|\le a_k^{-1}$ or $a_k\le |Th_k|\le b_k$
then $\|Th_k-\chi_GTh_k\|_Y \le (C+1)\eta.$ Let
$y_k=\chi_GTh_k.$ Then (6) and (7) follow.

This completes the inductive construction.  We now observe
that for every $2^N\le k\le 2^{N+1}-1$ we have
$(1-\gamma)^N\le \|x\chi_{F_k}\| \le \gamma^N.$ In
particular $\|h_k\|_X \ge (1-\gamma)^N$ for $0\le k\le
2^N-1.$ Thus $\|Th_k\|_Y \le (1-\gamma)^N\delta.$ By choice
of $\eta$ this implies that $ \frac12(1-\gamma)^N\delta\le
\|y_k\|_Y\le 1.$ Now we can appeal to Lemma 6.3 to deduce
that $(y_k)_{k=0}^{2^N-1}$ is 12-equivalent to a disjoint
sequence in $E_Y.$ In particular it is 12-unconditional.
Since $\|Th_k-y_k\|\|y_k\|^{-1} \le
2(C+1)\eta(1-\gamma)^{-N}\delta^{-1}$ we have $$
\sum_{k=0}^{2^N}\|Th_k-y_k\|\|y_k\|^{-1} \le
2^{N+1}(C+1)(1-\gamma)^{-N}\delta^{-1}\eta< 10^{-2}.$$

Hence $(Th_k)_{k=0}^{2^N-1}$ is 24-equivalent to a disjoint
sequence in $E_Y$ and hence $(h_k)_{k=0}^{2^N-1}$ is
$24\delta^{-1}$-equivalent to a disjoint sequence in $E_Y.$

We can define a linear map $L_N:CS_N\to X$ by
$L_N(\chi_{\Delta(\epsilon_1,\ldots\epsilon_N)})=
x\chi_{F_k}$ where $k= 2^N
+\frac12\sum_{j=1}^N(1-\epsilon_j)2^{N-j}.$ Then we can
induce a lattice norm on $CS_N$ by $\|f\|_N= \|L_Nf\|_X.$
Let $\Cal U$ be a non-principal ultrafilter on $\bold N.$ We
define for $f\in CS,$ $$ \|f\|_Z =\lim_{\Cal U}\|f\|_N.$$
Then $\|\,\|_Z$ is a lattice norm on $CS$ with the property
that if $E\in\Cal A_N$ then $(1-\gamma)^N\le \|\chi_E\|_Z
\le \gamma^N.$ Thus the completion $Z$ of this space is a
nonatomic Banach lattice which is finitely representable in
$X.$ Also the Haar system is clearly an unconditional basis
of $Z$ which is $25\delta^{-1}-$lattice finitely
representable in $E_Y.$\bull\enddemo

Before proving the next theorem, which is the main result of
the section, we make some definitions.  Let us denote by
$[0,\infty]$ the one-point compactification of $[0,\infty).$
Suppose $(\Omega_n,\mu_n)_{n=0}^{\infty}$ is a sequence of
Polish spaces with associated $\sigma-$finite measures and
let $f_n:\Omega_n\to [0,\infty]$ be Borel functions such
that for each $a>0$ we have $\mu_n(f_n>a)<\infty.$ We will
say that $(f_n,\mu_n)_{n=1}^{\infty}$ converges to
$(f_0,\mu_0)$ {\it in law} if and only if for every
continuous function $\phi:[0,\infty]\to \bold R$ so that
$\phi$ vanishes on a neighborhood of $0$ we have
$$\lim_{n\to\infty}\int_{\Omega_n}\phi\circ f_n\,d\mu_n=
\int_{\Omega_0}\phi\circ f_0\,d\mu_0.$$

If $f_n$ converges to $f_0$ in law then it is not difficult
to see that
$$\mu_0(f_0>a)\le\liminf_{n\to\infty}\mu_n(f_n>a)\le
\limsup_{n\to\infty}\mu_n(f_n>a) \le \mu_0(f_0\ge a).$$
Hence we can deduce that $f_n^*\to f_0^*$ a.e. on
$[0,\infty)$ and for any r.i. space $Y$ this implies that
$\|f_0\|_{Y(\Omega_0,\mu_0)}\le
\liminf_{n\to\infty}\|f_n\|_{Y(\Omega_n,\mu_n)}.$

\proclaim{Theorem 6.7} Suppose $Y$ is an r.i. space on
$[0,1]$ or $[0,\infty)$ with nontrivial concavity.  Suppose
$X$ is a good K\"othe function space on $(\Delta,\lambda)$
which is isomorphic to a subspace of $Y.$ Then
either:\newline (1) There is a nonatomic Banach lattice $Z$
which is lattice-finitely representable in $X$ and such that
$Z$ has an unconditional basis, which is lattice finitely
representable in $E_Y,$ or:\newline (2) There is a
cone-embedding of $X_{1/2}$ into $Y_{1/2}.$ \endproclaim

\demo{Proof}Let $\Cal C$ be the (countable) algebra of
clopen subsets of $\Delta.$ We define a compact space
$\Omega=[0,\infty]^{\Cal C}.$ We denote the co-ordinate maps
on $\Omega$ by $\xi_E$ for $E\in\Cal C.$

Let us suppose first that $Y=Y[0,\infty);$ we will describe
the minor modifications for the case $[0,1]$ afterwards.  We
suppose that $Y$ is $q$-concave with constant one where
$q<\infty.$ Suppose $p>2q$ is fixed.  Let $T:X\to Y$ be a
linear map satisfying for some $\delta>0,$ $\delta\|x\|_X\le
\|Tx\|_Y\le \|x\|_X$ for $x\in X,$ and define $Q_n:CS_n\to
L_0[0,\infty)$ as above.

We make first the observation that, as $Y$ is $q-$concave,
we have an estimate $\|\chi_{[0,t]}\|_Y \ge t^{1/q}$ for
$t\ge 1$ and hence if $y\in Y$ then $y^*(t)^q \le
t^{-1}\|y\|^q_Y$ for $t\ge 1.$ It follows that if $y\in Y$
then $$ \int_0^{\infty}\min(1,|y|^{p/2})\,dt \le
1+\|y\|_Y^{p/2q}\int_1^{\infty}t^{-p/2q}dt \le
1+C_0\|y\|_Y^{p/2q}$$ for a suitable constant
$C_0=C_0(q,p).$

Let us define $\kappa_n:[0,\infty)\to \Omega$ by
$\xi_E\circ\kappa_n = Q_n(\chi_E)$ if $E\in\Cal C_n$ and
$\xi_E\circ\kappa_n=0$ otherwise.  Let $w$ be the weight
function on $\Omega$ defined by $w=\min(1,\xi_{\Delta}^p).$
We will define a Borel measure $\nu_n$ on $\Omega$ by
$$\nu_n(B) =
\int_{\kappa_n^{-1}B}\min(1,Q_n(\chi_{\Delta})^p)d\lambda.$$

Let us first note that
$$\nu_n(\Omega)=\int\min(1,(Q_n\chi_{\Delta})^p)d\lambda \le
1+C_0K_G^{p/q}$$ so that the sequence of Borel measures
$(\nu_n)$ is bounded in $\Cal M(\Omega).$ It follows that
$(\nu_n)$ has a weak$^*$-limit point $\nu.$ Let us define
$\mu_n=w^{-1}\nu_n$ and $\mu=w^{-1}\nu$; these measures are
$\sigma-$finite.

Note first that if $U$ is an open subset of $\Omega$ then
$\nu(U)\le\limsup \nu_n(U).$ We use this first to argue that
$\xi_E<\infty,\ \mu-$a.e. for every $E\in\Cal C.$ In fact if
$a>0$ then in $E\in \Cal C_n, $ we have $\nu_n(\xi_E>a) \le
\lambda(Q_n(\chi_E)>a)$ and by Lemma 6.4,
$a^{1/2}\min(1,\lambda(Q_n(\chi_E)>a)) \le K_G\|\chi_E\|_X.$
Hence $\lim_{a\to\infty}\nu(\xi_E>a)=0$ and so
$\mu(\xi_E=\infty)=0.$

Next we argue that if $E,F\in \Cal C$ are disjoint then
$\xi_{E\cup F}=\xi_E+\xi_F$ a.e. for $\mu.$ In fact, if
$\varepsilon>0,$ let $U$ be the set of $\omega\in\Omega$
such that $\xi_E(\omega),\xi_F(\omega),\xi_{E\cup
F}(\omega)<\infty$ and
$|\xi_E(\omega)+\xi_F(\omega)-\xi_{E\cup
F}(\omega)|>\varepsilon.$ Then if $E,F\in\Cal C_n, $ we have
$\nu_n(U)=0.$ Hence $\nu(U)=0$ and $\mu(U)=0.$ Thus
$\xi_E+\xi_F=\xi_{E\cup F}$ a.e.  It follows that we can
define a linear map $S_0:CS\to L_0(\mu)$ by
$S_0(\chi_E)=\xi_E.$

Now suppose $f\in CS_+.$ Let
$f=\sum_{k=1}^N\alpha_k\chi_{E_k}$ where $E_1,\ldots,E_N$
are clopen sets in $\Delta,$ and $\alpha_k\ge 0$ for $1\le
k\le N.$ Let $g=\sum_{k=1}^N\alpha_k\xi_{E_k}$ so that
$g=S_0f$ a.e. for $\mu.$ Let $M=\sum_{k=1}^N\alpha_k.$ Then
$f\le M\chi_{\Delta}$ and $g\le M\xi_{\Delta},$ a.e. for
$\mu.$

For any $a>0$, let $\varphi_a$ be a continuous function on
$[0,\infty]$ such that $\varphi_a(t)=0$ if $0\le t\le
1/(2Ma)$ and $\varphi_a(t)=1$ if $t\ge 1/(Ma).$ Then let
$g_a=(\varphi_a\circ\xi_{\Delta})\min(a,g).$ Then $\tau_ag
\le g_a\le g,$ $\mu$-a.e.

For fixed $a>0,$ $g_a$ is continuous on $\Omega.$
Furthermore for each $n,$ $\mu_n(g_a>0) \le
\mu_n(\xi_{\Delta}>(Ma)^{-1})\le
\lambda(Q_n(\chi_{\Delta})>(Ma)^{-1})$ is uniformly bounded.
If $\nu_{n(k)}$ converges weak$^*$ to $\nu$ then for any
continuous function $\varphi$ on $[0,\infty]$ which vanishes
in a neighborhood of the origin, we have $$ \align
\lim_{k\to\infty}\int_{\Omega}g_a\,d\mu_{n_k} &=
\lim_{k\to\infty}\int_{\Omega}g_aw^{-1}d\nu_{n_k}\\
&=\lim_{k\to\infty}\int_{\Omega}
(\varphi_a\circ\xi_{\Delta})
\max(1,\xi_{\Delta}^{-p})\min(g,a)d\nu_{n_k} \\
&=\int_{\Omega} (\varphi_a\circ\xi_{\Delta})
\max(1,\xi_{\Delta}^{-p})\min(g,a)d\nu \\ &=
\int_{\Omega}g_ad\mu.  \endalign $$ Thus $(g_a,\mu_{n_k})$
converges in law to $(g_a,\mu).$ Since $Y$ is
order-continuous, $g_a$ is bounded and the measures of the
supports are uniformly bounded, this implies that
$$\lim_{k\to\infty}\|g_a^{1/2}\|_{Y(\mu_{n_k})}=\|g_a^{1/2}\|_{Y(\mu)}.$$

If $E_1,\ldots,E_N\in\Cal C_n$ then we have $g\ge g_a\ge
\tau_ag$ a.e. for $\mu_n.$ It follows that we have
$\|\tau_ag\|_{Y(\mu_n)} \le \|g_a\|_{Y(\mu_n)}\le
\|g\|_{Y(\mu_n)}.$

Note however that $(g,\mu_n)$ coincides in law with
$(Q_nf,\lambda)$ for if $B$ is a Borel subset of
$(0,\infty)$ then
$\mu_n(g^{-1}B)=\int_{g^{-1}B}w^{-1}d\nu_n=\lambda(\kappa_n^{-1}g^{-1}B)
=\lambda((Q_nf)^{-1}B).$

Hence we obtain the estimate
$$\lim_{a\to\infty}\liminf_{n\to\infty}\|\tau_a(Q_nf)^{1/2}\|_Y
\le \|g^{1/2}\|_{Y(\mu)}\le
\limsup_{n\to\infty}\|(Q_nf)^{1/2}\|_Y.$$

We conclude that $f\in CS_+$ we have $\Psi(f^{1/2})^2 \le
\|S_0f\|_{Y_{1/2}(\mu)}\le K_G^2\|f\|_{X_{1/2}}.$ Thus $S_0$
extends to a bounded positive operator $S:X_{1/2}\to
Y_{1/2}.$ If alternative (1) of the theorem is false then,
by Lemma 6.6, $S$ has a lower estimate and it is clear that
$S$ is a cone-embedding, as required.

In the case when $Y=Y[0,1]$ we can regard $Y$ as being
embedded in a space modelled on $[0,\infty)$ and need only
observe that in the above proof, the measures $\mu_n$ and
$\mu$ have total mass at most one.  \bull \enddemo

\proclaim{Theorem 6.8} Suppose $Y$ is an r.i. space on
$[0,1]$ or $[0,\infty)$ with nontrivial concavity, which is
either strictly 2-convex or of Orlicz-Lorentz type.  Suppose
$X$ is a good K\"othe function space on $(\Delta,\lambda)$
which is isomorphic to a subspace of $Y.$ Then there is a
cone-embedding of $X_{1/2}$ into $Y_{1/2}.$ \endproclaim

\demo{Proof}It is enough to show that the existence of $Z$
in Theorem 6.7 leads to a contradiction.  Suppose first that
$Y$ is strictly 2-convex; then $E_Y$ is also strictly
2-convex.  This implies that the unconditional basis of $Z$
is strictly 2-convex, and hence $Z$ can contain no copy of
$\ell_2$; however $Z$ must have nontrivial cotype and this
contradicts Lemma 2.4 of \cite{11}.

If $Y$ is of Orlicz-Lorentz type then $E_Y$ is
lattice-isomorphic to a modular sequence space which has
nontrivial cotype.  Now the unconditional basis of $Z$ is
lattice finitely representable in $E_Y.$ This implies that
$Z$ also is isomorphic to a modular sequence space, also
with nontrivial cotype.  This can be established directly
without difficulty, but is also a special case of more
general results on ultraproducts of Orlicz spaces and
Orlicz-Musielak spaces, for which we refer to \cite{12},
\cite{18} and \cite{43}.  This now contradicts Theorem 4.3
and Corollary 4.4 in \cite{28} (which in turn extends an
earlier result of Lindenstrauss and Tzafriri
\cite{32}).\bull \enddemo

\vskip10pt

\heading{7.  The main results}\endheading \vskip10pt

Before proving our main results for embeddings of nonatomic
Banach lattices into r.i. spaces, we first give an
illustrative theorem for atomic Banach lattices.  Compare
this result with those of Johnson and Schechtman \cite{22}
and Carothers and Dilworth \cite{8}.

\proclaim{Theorem 7.1}Suppose $Y$ is an r.i. space on
$[0,\infty)$ with nontrivial cotype, and suppose that either
(a) $Y$ is 2-convex or (b) $p_Y>2.$ Suppose $(u_n)$ is a
strictly 2-convex unconditional basic sequence in $Y.$ Then
$(u_n)$ is equivalent to a disjoint sequence.  Equivalently,
if $X$ is a strictly 2-convex atomic Banach lattice which is
isomorphic to a subspace of $Y$ then $X$ is
lattice-isomorphic to a sublattice of $Y.$\endproclaim

\demo{Remark}We do not know if this theorem holds when $Y$
is an r.i. space on $[0,1].$\bull\enddemo

\demo{Proof}Let us suppose that $X$ is an atomic Banach
lattice represented as a function space of $\bold N$ with
canonical basis vectors $e_n$ and that $S:X\to Y$ is an
embedding with $Se_n=u_n.$ Then by Theorem 1.d.6 of
\cite{34} we can define a cone-embedding $L:X_{1/2}\to
Y_{1/2}$ by $Le_n=|u_n|^2.$ The result is now obtained by
putting together the facts previously established on
cone-embeddings.  Since $X_{1/2}$ is strictly 1-convex $L$
is a strong cone-embedding, by Lemma 4.1; then since
$Y_{1/2}$ has property (d) under either conditions (a) or
(b), Proposition 4.5 shows that $X_{1/2}$ is
lattice-isomorphic to a sublattice of $Y_{1/2}.$ But this
implies the result.\bull\enddemo

We now prove the nonatomic version of the above theorem.

\proclaim{Theorem 7.2}Suppose $Y$ be an r.i. space on
$[0,\infty)$ with nontrivial concavity and either \newline
(a) $Y$ is strictly 2-convex, or\newline (b) $Y$ is 2-convex
and of Orlicz-Lorentz type, or \newline (c) $p_Y>2$ and $Y$
is of Orlicz-Lorentz type.  Suppose $X$ be a strictly
2-convex nonatomic Banach lattice.  If $X$ is isomorphic to
a subspace of $Y$, then $X$ is isomorphic to a sublattice of
$Y.$ \endproclaim

\demo{Proof}We can of course assume that $X$ is a good
K\"othe function space on $(\Delta,\lambda)$.  We first
apply Theorem 6.8 to deduce the existence of a
cone-embedding of $X_{1/2}$ into $Y_{1/2}.$ Now the proof
proceeds as in Theorem 7.1.  \bull\enddemo

\demo{Remarks} Let us first note that if $Y$ is 2-convex
then $X$ must also be 2-convex at least; the hypothesis that
$X$ is strictly 2-convex is then equivalent to the
hypothesis that $\ell_2$ is not lattice finitely
representable in $X$ (cf.  \cite{21} Lemma 2.4).  This
result was previously known in the special case
$Y=L_p[0,\infty)$ \cite{21} ,Theorem 1.8 (the atomic case is
proved in \cite{16}.)  \enddemo

We now turn to the case when $Y$ is an r.i. space on
$[0,1];$ here our result is not quite as strong (exactly as
in the atomic case:  see discussion after Theorem 7.1).

\proclaim{Theorem 7.3}Let $Y$ be an r.i. space on $[0,1]$
with nontrivial concavity and suppose either (a) $Y$ is
strictly 2-convex or (b) $p_Y>2$ and $Y$ is of
Orlicz-Lorentz type.  Suppose $X$ is a nonatomic strictly
2-convex Banach lattice which is isomorphic to a subspace of
$Y$.  Then $X$ contains a nontrivial band $X_0$ which is
lattice-isomorphic to a sublattice of $Y.$ \endproclaim

\demo{Proof}We will consider $X$ as a good K\"othe function
space on $[0,1].$ Then there is, by Theorem 6.8, a
cone-embedding $L:X_{1/2}\to Y_{1/2}.$ Furthermore $X_{1/2}$
is $s$-convex for some $s>1$ and there exists in either case
$r>2$ so that $Y_{1/r}$ has property (d).  Proposition 4.6
then implies that for some Borel set $E$ with $\lambda(E)>0$
the band $X_{1/2}(E)$ is lattice-isomorphic to a sublattice
of $Y_{1/2}.$ The result then follows.\bull\enddemo

We now turn our attention to the case when $X$ is known to
be an r.i. space.

\proclaim{Corollary 7.4} Let $Y$ be an r.i. space on
$I=[0,1]$ or $[0,\infty)$ with nontrivial concavity.
Suppose either \newline (a) $Y$ is strictly 2-convex or
\newline (b) $Y$ is of Orlicz-Lorentz type and $p_Y>2.$
\newline Suppose $X$ is an r.i. space on $I=[0,1],$ with
$X\neq L_2[0,1].$ Assume that $X$ is isomorphic to a
subspace of $Y.$ Then $X$ is isomorphic to a sublattice of
$Y$ and there exists $f\in Y$ so that $X=Y_f[0,1].$
\endproclaim

\demo{Proof}Consider first case (a).  By Proposition 2.e.10
of \cite{28} or Section 2 of \cite{21} either $X=L_2$ or $X$
is strictly 2-convex.  The result then follows by the
preceding Theorems 7.2 and 7.3.

Case (b) is slightly different.  In this case Theorem 6.8
implies that there is a cone-embedding of $X_{1/2}$ into
$Y_{1/2}.$ By Proposition 5.2, either $X_{1/2}=L_1$ (i.e.
$X=L_2$) or $X_{1/2}$ is isomorphic to a sublattice of
$Y_{1/2}$ and the result follows.\bull\enddemo

\demo{Remarks}Some special cases of Corollary 7.4 are known.
In \cite{21} Theorem 7.7 the corollary is proved when $Y$ is
a strictly 2-convex Orlicz function space.  Later, Carothers
\cite{5} and \cite{6} proves the same theorem for Lorentz
spaces $L_{p,q}$ where $p>\max(q,2)$.  Carothers considers
first the strictly 2-convex case ($2\le q< p$) and later
modifies the proof to the case $1\le q\le 2 <p.$ Note that
in these cases and in more general Lorentz spaces considered
by Carothers one has the additional information that every
$Y_f[0,1]$ coincides with $Y[0,1].$ This is equivalent to an
inequality of the form $\|f\otimes g\|_Y \le
K\|f\|_Y\|g\|_Y$ for $f,g\in Y[0,1].$ This additional
information is actually used in the proof.  \enddemo

For reference let us state one additional case which follows
from Theorem 7.2 and Proposition 3.3.

\proclaim{Corollary 7.5}Let $Y$ be an r.i. space on
$[0,\infty)$ with nontrivial concavity which is 2-convex and
of Orlicz-Lorentz type.  Let $X$ be a strictly 2-convex r.i.
space on $[0,1]$ which is isomorphic to a subspace of $Y.$
Then there exists $f\in Y$ so that $X=Y_f[0,1].$\endproclaim

Let us note the following special case.

\proclaim{Corollary 7.6} Suppose $2<p<\infty$ and $Y$ is a
$p$-convex r.i. space on $[0,1]$ or $[0,\infty)$ with
nontrivial concavity.  Suppose $L_p$ is isomorphic to a
subspace of $Y.$ Then $Y[0,1]=L_p[0,1].$\endproclaim

\demo{Proof}It follows from Corollary 7.3 that
$L_p[0,1]=Y_f[0,1] \subset Y[0,1]$ but $Y[0,1]\subset
L_p[0,1]$ since $Y$ is $p$-convex.\bull\enddemo

\demo{Remarks} The condition that $Y$ is $p$-convex cannot be relaxed
here (cf. \cite{19}).  We remark that analogues of Corollary
7.6 for
$1\le p<2$ have been proved in several places in the
literature.  In the case $p=1,$ then $L_1$ embeds into a
separable r.i. space $Y[0,1]$ if and only if
$Y[0,1]=L_1[0,1].$ This is proved under the additional
hypothesis that $Y$ has nontrivial cotype in \cite{21} (cf.
\cite{34} Corollary 2.e.4); it is proved under the
hypothesis that $Y$ does not contain $c_0$ in \cite{23}.
The result with no additional hypothesis follows from
Theorem 10.7 and Theorem 7.3 of \cite{27}.  For the case
$1<p<2$ a similar result holds when $Y$ is separable and
$p$-convex provided one eliminates the possibility that $Y$
contains a disjoint sequence equivalent to the Haar basis of
$L_p[0,1]$ (see Theorems 7.3 and 10.7 of \cite{27
}.)\enddemo

In our final result we consider the case when instead $Y$ is
$p$-concave for some $p>2$ and $L_p$ embeds into $X$.

\proclaim{Theorem 7.7}Suppose $2<p<\infty$ and that $Y$ is a
$p$-concave r.i. space on $[0,1]$ or $[0,\infty).$ Suppose
that $L_p$ is isomorphic to a subspace of $Y.$ Then,
either:\newline (a) The Haar basis of $L_p$ is lattice
finitely-representable in $E_Y$ or\newline (b)
$Y[0,1]=L_p[0,1].$

In particular, if $Y$ is strictly 2-convex or of
Orlicz-Lorentz type, then $Y[0,1]=L_p[0,1].$ \endproclaim

\demo{Proof} We will apply Theorem 6.7.  First suppose that
$Z$ is a nonatomic Banach lattice which is lattice finitely
representable in $L_p,$ which has an unconditional basis
lattice finitely representable in $E_Y.$ Then of course
$Z=L_p$.  It follows from the reproducibility of the Haar
basis (Theorem 2.c.8 of \cite{34}) that the Haar basis is
also lattice finitely representable in $E_Y,$ contrary to
hypothesis.

We conclude that $L_{p/2}$ can be cone-embedded into
$Y_{1/2}.$ Now the result follows immediately from
Proposition 5.3.\bull\enddemo

\demo{Remarks}Here, the condition that $Y$ is $p$-concave cannot
be relaxed (\cite{19}).   We give a simple application.  Suppose
$1\le r<2<p$ and
$Y=(L_r+L_p)[0,\infty).$ It follows from the above theorem
that $L_p$ is not isomorphic to a subspace of $Y$ which
answers a question raised in \cite{17}.\enddemo

\vskip10pt \heading{8.  Complemented subspaces of r.i.
spaces}\endheading \vskip10pt

The following result is quickly deduced from the methods of
\cite{27}.

\proclaim{Theorem 8.1}Let $Y$ be a separable
order-continuous Banach lattice, which contains no
complemented sublattice isomorphic to $\ell_2.$ Suppose $X$
is a Banach lattice which is isomorphic to a complemented
subspace of $Y$.  Then either:\newline (a) There is a
constant $C$ so that, for every $n,$ $\ell_2^n$ is
C-lattice-isomorphic to a complemented sublattice of $X$,
or:\newline (b) There exists $N$ so that $X$ is lattice
isomorphic to a complemented sublattice of
$Y^N=Y\oplus\cdots\oplus Y.$ \endproclaim

\demo{Proof}We will prove under the assumption that $X$ is
nonatomic.  (An exposition of the atomic case, which is
proved by the same techniques, will be given in \cite{10}.)
In this case we may suppose that both $X$ and $Y$ are good
K\"othe function spaces on $(\Delta,\lambda)$ and that $X$
has the ``strong density property.''  By combining Theorems
6.1 and 6.3 of \cite{27} it is possible to find a sequence
of Borel maps $\sigma_n:\Delta\to \Delta$ and three
sequences $(a_n^P),(a_n^Q),(a_n^R)$ of nonnegative Borel
functions on $\Delta$ so that $a_n^P(s)^2\le
a_n^Q(s)a_n^R(s)$ and if:  $$ \align Pf
&=\sum_{n=1}^{\infty}a_n^P f\circ\sigma_n\\ Qf
&=\sum_{n=1}^{\infty}a_n^Q f\circ\sigma_n\\ Rf
&=\sum_{n=1}^{\infty}a_n^R f\circ\sigma_n \endalign $$ for
$f\in (L_0)_+$ then we have for a suitable constant $C_1$
that $\|Pf\|_1\le C_1\|f\|_1,$ $\|Qf\|_{Y_{1/2}}\le
C_1\|f\|_{X_{1/2}},$ and $\|Rf\|_{Y^*_{1/2}} \le
C_1\|f\|_{X^*_{1/2}}.$ Note here that $Q$ need only map into
$Y_{max,1/2}$ and not necessarily into $Y_{1/2}.$ Now by
Theorem 6.4 of \cite{27} it can be seen that if the first
alternative fails then there is a constant $c>0$ so that $$
\int\sup_n a_n^Pf\circ\sigma_nd\lambda \ge c\int f
d\lambda$$ for $f\ge 0.$ We now use an argument due to Dor
\cite{15}.  Consider the map $T:L_1\to L_1(c_0)$ defined by
$Tf(s) = (a_n^P(s)f(\sigma_n(s)).$ Then $\|T\|\le C_1$ and
$\|Tf\|\ge c\|f\|.$ Note that since $c_0$ has separable
dual, $L_1(c_0)^*$ can be identified with
$L_{\infty}(\ell_1)$.  By the Hahn Banach theorem there
exist $\phi_n\in L_{\infty}$ so that
$\|\sum_{n=1}^{\infty}|\phi_n|\|_{\infty}\le C_1c^{-1}$ and
$$ \sum_{n=1}^{\infty}\int \phi_na_n^Pf\circ\sigma_n
d\lambda =\int f\,d\lambda$$ for $f\in L_1(\lambda).$

Now for each $n$ define $E_n=\{s:\phi_n(s)>(2C_1)^{-1}\}.$
Then for $f\ge 0,$ $$
\sum_{n=1}^{\infty}\int_{\Delta\setminus
E_n}\phi_na_n^Pf\circ\sigma_n\le \frac12 \int f\,d\lambda.$$
Hence $$
\sum_{n=1}^{\infty}\int_{E_n}a_n^Pf\circ\sigma_nd\lambda \ge
\frac{c}{2C_1}\int f\,d\lambda.$$

Notice that $\sum_{n=1}^{\infty}\chi_{E_n} \le
2C_1\sum_{n=1}^{\infty} |\phi_n|\le 2C_1^2c^{-1}$ almost
everywhere.  Let $N$ be the least integer greater than
$2C_1^2c^{-1}.$ Consider the operators $P'$, $Q'$ and $R'$
defined by $$ \align P'f
&=\sum_{n=1}^{\infty}a_n^P\chi_{E_n}f\circ\sigma_n\\ Q'f
&=\sum_{n=1}^{\infty}a_n^Q\chi_{E_n}f\circ\sigma_n\\ R'f
&=\sum_{n=1}^{\infty}a_n^R\chi_{E_n}f\circ\sigma_n.
\endalign $$ Then these operators can each be rewritten in
the form, $$ \align P'f &= \sum_{n=1}^{N}b_n^Pf\circ\pi_n\\
Q'f &= \sum_{n=1}^{N}b_n^Qf\circ\pi_n\\ R'f &=
\sum_{n=1}^{N}b_n^Rf\circ\pi_n.  \endalign $$ for suitable
nonnegative Borel functions $b_n^P,b_n^Q,b_n^R,$ for $1\le
n\le N,$ which also satisfy $(b_n^P)^2 \le b_n^Qb_n^R$ a.e.,
and for suitable Borel maps $\pi_n:\Delta\to\Delta.$

Now define $U:X\to Y_{max}^N,$ $V:X^*\to (Y^*)^N,$ by $$
\align Uf(s,n) &= (b_n^Q(s))^{1/2}f(\pi_n(s))\\ Vf(s,n) &=
(b_n^R(s))^{1/2}f(\pi_n(s)).  \endalign $$ It is easy to see
that $U$ is bounded for $$ \align \max_{1\le n\le
N}\|Uf(.,n)\|_Y &\le \|(\sum_{n=1}^N b_n^Q
(f\circ\pi_n)^2)^{1/2}\|_Y \\ &= \|Q'f^2\|_{Y_{1/2}}^{1/2}\\
&\le \|Qf^2\|_{Y_{1/2}}^{1/2}\\ &\le C_1^{1/2}\|f\|_X.
\endalign $$ Similarly $V$ is bounded.

The proof is completed by Proposition 2.3 of \cite{27}, for
if $F$ is a Borel subset of $\Delta$ then $$ \align
\sum_{n=1}^N \int_{\sigma_n^{-1}F}b_n^Q(s)^{1/2}
b_n^R(s)^{1/2}d\lambda &\ge \sum_{n=1}^N
\int_{\sigma_n^{-1}F}b_n^P(s)\, d\lambda\\ &=
\int_{\Delta}P'\chi_F(s)d\lambda\\ &\ge
\frac{c}{2C_1}\lambda(F).\bull \endalign $$ \enddemo

This theorem has immediate consequences if $Y$ is an r.i.
space.

\proclaim{Theorem 8.2}Let $Y$ be a separable r.i. space on
$[0,1]$ or $[0,\infty),$ which contains no complemented
sublattice isomorphic to $\ell_2.$ Suppose $X$ is a strictly
2-convex or strictly 2-concave Banach lattice which is
isomorphic to a complemented subspace of $Y.$ Then $X$ is
lattice-isomorphic to a complemented sublattice of $Y.$
\endproclaim

We remark that Theorem 8.2 is closely related to Theorem 8.1
of \cite{27}, and could be used to simplify some of the
arguments in the proof of that theorem somewhat.

\enddocument
\bye